\input amstex 
\documentstyle{amsppt}
\loadbold
\magnification=\magstep1
\pageheight{9.0truein}
\pagewidth{6.5truein}
\NoBlackBoxes
\TagsAsMath

\def\ZZ{{\Bbb Z}}

\def\CC{{\Bbb C}}
\def\AA{{\Bbb A}}
\def\B{{\Cal B}}

\def\P{{\Cal P}}

\def\U{{\Cal U}}
\def\SS{{\Cal S}}
\def\b{{\frak b}}

\def\e{{\frak e}}
\def\f{{\frak f}}
\def\g{{\frak g}}
\def\h{{\frak h}}
\def\p{{\frak p}}
\def\l{{\frak l}}
\def\n{{\frak n}}

\def\so{{\frak{so}}}
\def\sll{{\frak{sl}}}
\def\sp{{\frak{sp}}}

\def\wm{{w_\circ}}             
\def\wJ{{w^J_\circ}}
\def\dw{\dot{w}}

\def\Ga{\Gamma}                
\def\De{\Delta}
\def\al{\alpha}

\def\be{\beta}
\def\ga{\gamma}
\def\Om{\Omega}
\def \WJmax {W^J_{\max}}
\def \WJmin {W^J_{\min}}
\def \ra{\rightarrow}             
\def \hra{\hookrightarrow}
\def \mt{\mapsto}

\def\wt{\widetilde}
\def\ol{\overline}

\def\leftmatrix{\left[ \smallmatrix \format \c&&\quad\c \\}
\def\rightmatrix{\endsmallmatrix\right]}

\def\CCx{\CC^{\times}}

\def\ci{\circ}

\def\Lie{\operatorname{Lie}}
\def\Cl{\operatorname{Cl}}
\def\supp{\operatorname{supp}}

\def\Span{\operatorname{Span}}
\def\sign{\operatorname{sign}}

\def\Ad{\operatorname{Ad}}

\def\ad{\operatorname{ad}}
\def\id{\operatorname{id}}

\def\diag{\operatorname{diag}}

\def\Bour{{\bf 1}}
\def\BL{{\bf 2}}
\def\BGY{{\bf 3}}
\def\Car{{\bf 4}}
\def\CF{{\bf 5}}
\def\Deo{{\bf 6}}
\def\Dri{{\bf 7}}
\def\FL{{\bf 8}}
\def\FH{{\bf 9}}
\def\G{{\bf 10}}
\def\Kam{{\bf 11}}
\def\LY{{\bf 12}}
\def\Lus{{\bf 13}}
\def\MRS{{\bf 14}}
\def\Mus{{\bf 15}}
\def\Nou{{\bf 16}}
\def\RTF{{\bf 17}}
\def\RRS{{\bf 18}}
\def\Rie{{\bf 19}}
\def\Sp{{\bf 20}}
\def\St{{\bf 21}}
\def\Ta{{\bf 22}}
\def\Wi{{\bf 23}}
\def\Wone{{\bf 24}}
\def\Wtwo{{\bf 25}}
\def\X{{\bf 26}}
\topmatter

\title Poisson structures on affine spaces and flag varieties. II.
General case
\endtitle

\rightheadtext{Poisson structures on flag varieties}

\dedicatory Dedicated to the memory of our colleague Xu-Dong Liu (1962-2005)
\enddedicatory

\author K. R. Goodearl and M. Yakimov
\endauthor

\address Department of Mathematics, University of California,
Santa Barbara, CA 93106, USA\endaddress
\email goodearl\@math.ucsb.edu\endemail

\address Department of Mathematics, University of California,
Santa Barbara, CA 93106, USA\endaddress
\email yakimov\@math.ucsb.edu\endemail


\thanks 
The research of the first author was partially supported by 
National Science Foundation grant DMS-0401558.
The research of the second author was partially supported 
by National Science Foundation grant DMS-0406057 and
an Alfred P. Sloan research fellowship. 
\endthanks

\abstract The standard Poisson structures on the flag varieties 
$G/P$ of a complex reductive algebraic group $G$ are investigated. 
It is shown that the orbits of symplectic leaves in $G/P$ 
under a fixed maximal torus of $G$ are smooth irreducible
locally closed subvarieties of $G/P$, isomorphic to intersections
of dual Schubert cells in the full flag variety $G/B$ of $G$,
and their Zariski closures are explicitly computed. Two different
proofs of the former result are presented. The first 
is in the framework of Poisson homogeneous spaces and the 
second one uses an idea of weak splittings of surjective
Poisson submersions, based on the notion of Poisson--Dirac 
submanifolds. For a parabolic subgroup $P$ 
with abelian unipotent radical (in which case $G/P$ is
a Hermitian symmetric space of compact type), it is 
shown that all orbits of the standard Levi factor $L$ 
of $P$ on $G/P$ are complete Poisson subvarieties which 
are quotients of $L$, equipped with the standard 
Poisson structure. Moreover, it is proved that the 
Poisson structure on $G/P$ vanishes at all special
base points for the $L$-orbits on $G/P$ constructed
by Richardson, R\"ohrle, and Steinberg.  
\endabstract

\endtopmatter
\document


\head Introduction \endhead
\definition{0.1}
First we fix some notation. Let $G$ be a connected reductive algebraic 
group over $\CC.$ Fix a pair of dual Borel subgroups $B^\pm$ in 
$G$, and set $H = B^+ \cap B^-$ for the corresponding maximal torus.
Denote the corresponding Lie algebras by $\g$, $\b^\pm$, and $\h$.
Denote by $\De$ and $\De^+$ the set of all roots, 
respectively 
all positive roots, of $\g$ with respect to $\h$. 
Let $\Ga$ be the set of all positive simple roots of $\g$. 
For a subset $J$ of $\Ga$, let $P_J^\pm$ 
be the standard parabolic subgroups of $G$, containing respectively 
the Borel subgroups $B^\pm$. Let $L_J=P_J^+ \cap P_J^-$ be the common
Levi factor of
$P_J^\pm$. Denote by $U^\pm$ and $U_J^\pm$  
the unipotent radicals of $B^\pm$ and $P_J^\pm,$ respectively. 
Set $\p_J^\pm = \Lie P_J^\pm$, $\l_J = \Lie L_J$, 
$\n^\pm = \Lie U^\pm$, 
and $\n_J^\pm = \Lie U_J^\pm$. 
\enddefinition

\definition{0.2} We fix a nondegenerate bilinear invariant form 
$\langle.,.\rangle$ on $\g$ for which the square of the length of 
a long root is equal to 2. 
Recall \cite{\BGY, eq\. (1.1)} that the 
standard $r$-matrix of $\g$ is given by
$$
r_\g = \sum_\al \frac{\langle \al, \al \rangle}{2} 
e_\al \wedge f_\al 
\tag 0.1
$$
where
$\{e_\al\}_{\al \in \De^+}$ and $\{f_\al\}_{\al \in \De^+}$
are any sets of root vectors of $\g$, 
normalized by 
$[e_\al, f_\al] = \al\spcheck= 
2 \al /\langle \al, \al \rangle$.
(In the last equation, $\h$ and $\h^*$ are identified 
via the restriction of the form $\langle.,.\rangle$;
for $\al \in \De^+$, $e_\al$ and $f_\al$ are root vectors 
for the roots $\al$ and $-\al$, respectively.)
The standard Poisson structure on $G$ is given by
$$
\pi_G = L(r_\g) - R(r_\g)= \chi^R(r_\g) - \chi^L(r_\g).
\tag 0.2
$$ 
Here $L(r_\g)$ and $R(r_\g)$ refer to the left and right 
invariant bivector fields on $G$ associated 
to $r_\g \in \wedge^2 \g \cong \wedge^2 T_e G$. 

For any subset $J \subseteq \Ga$, the standard parabolic subgroup $P_J^+$
is a Poisson algebraic subgroup of $(G, \pi_\g)$. The natural
projection 
$$
\eta_J : G \ra G/P_J^+
\tag 0.3
$$ 
induces the following Poisson
structure on the flag variety $G/P_J^+$: 
$$
\pi_J := \eta_{J *}(\pi) = - \chi(r_\g),
\tag 0.4
$$ 
see \cite{\BGY, Theorem 1.8} and Proposition 1.3 below.
Throughout the paper, for an action of $G$ 
on a variety $M$, we denote the extension to 
$\wedge \g$ of the infinitesimal action of $\g$ on $M$ 
by
$$
\chi : \wedge \g \ra \Ga(M, \wedge TM).
$$ 
For brevity, we set $\eta = \eta_\emptyset$ and 
$\pi=\pi_\emptyset$.

Since the Poisson structure $\pi_G$ vanishes on $H$, 
the left action of $H$ on $G/P_J^+$ preserves
$\pi_J$.
\enddefinition

\definition{0.3} In this paper we investigate the geometry of the
Poisson structures $\pi_J$ on the flag varieties $G/P_J^+$. This
continues our work with Brown in \cite{\BGY}, which we shall refer to as
Part I.

Before we state the main results of the paper, we 
introduce some notation on Weyl groups. 
The Weyl group of 
the pair $(G,H)$ will be denoted by $W$. 
For $w \in W$, we will denote  by $\dw$ a representative 
of $w$ in the normalizer of $H$. If a formula 
does not depend on the choice of this representative
the dot will be omitted. The Weyl group
of $(L_J,H)$, naturally thought of as a subgroup of $W$, 
will be denoted by $W_J$. Recall that 
each coset in $W/W_J$ has unique minimal and maximal length 
representatives. The sets of those will be denoted by $\WJmin$ and
$\WJmax$ respectively. Denote the longest elements 
of $W$ and $W_J$ by $\wm$ and $\wJ$. 

For $w \in W$ and $J \subseteq \Ga$, set 
$$
x_w^J = w P_J^+ \in G/P_J^+, \quad \qquad
x_w = x_w^\emptyset = w B^+ \in G/B^+.
\tag 0.5
$$
\enddefinition

\definition{0.4} The following Theorem summarizes some of
our results (Theorems 1.5, 1.8, 4.6, and Proposition 4.2).
Recall from Part I, \cite{\BGY, eq\. (2.11)},
the notation
$$
\U_{\dw_1, w_2} = U_{w_1}^- \dw_1 \cap B^+ w_2 B^+, 
\quad w_1, w_2 \in W
\tag 0.6 
$$
where for $w \in W$,
$$
U_w^- = U^- \cap \Ad_w(U^-).
\tag 0.7
$$
\enddefinition

\proclaim{Theorem} {\rm{(i)}} There are only finitely many $H$-orbits 
of symplectic leaves on $(G/P_J^+, \pi_J)$, parametrized
by pairs $(w_1, w_2) \in \WJmax \times W$ such that 
$w_1 \leq w_2$ in the Bruhat order. The torus orbit
corresponding to the pair $(w_1, w_2)$ is given by
$$
\SS_{w_1,w_2}^J = \U_{\dw_1, w_2} . P_J^+,
$$  
and is biregularly isomorphic to 
the intersection $\B_{w_1,w_2}= B^- x_{w_1} \cap B^+ x_{w_2}$
of dual Schubert cells in the generalized
full flag variety $G/B^+$.
Thus, the $H$-orbits of symplectic leaves
on $(G/B^+, \pi)$ are exactly the 
intersections of dual Schubert cells
$\B_{w_1,w_2}$.

{\rm{(ii)}} For each pair $(w_1, w_2)$ as above,
the Zariski closure of $\SS^J_{w_1, w_2}$
in $(G/P_J^+, \pi_J)$ is equal to the 
union of $\SS^J_{v_1, v_2}$ over those pairs
$(v_1, v_2) \in \WJmax \times W$ with $v_1 \leq v_2$
for which there exists
$z \in W_J$ such that 
$w_1 \leq v_1 z$ and 
$w_2 \geq v_2 z$.   

{\rm{(iii)}} If $P_J^+$ is a parabolic subgroup of $G$ with abelian
unipotent radical {\rm(}in which case $G/P_J^+$ is a hermitian
symmetric space of compact type{\rm)}, then all $L_J$-orbits
on $G/P_J^+$ are complete Poisson subvarieties of 
$(G/P_J^+, \pi_J)$ and all of them are quotients 
of the standard Poisson group structure on $L_J$. 
\endproclaim

In addition, we prove that in the case of Hermitian symmetric 
spaces of compact type the Poisson structure $\pi_J$ vanishes 
at the base points for the $L_J$-orbits constructed by 
Richardson, R\"ohrle, and Steinberg \cite{\RRS}. 
In this case, we also characterize explicitly the 
$H$-orbits of symplectic leaves which fall within 
a given $L_J$-orbit. This relates part (i) and (iii) 
of the above Theorem. It is done in Theorem 4.13 and 
will not be formulated here due to the amount of 
notation it requires.

The partition of the partial flag variety $G/P_J^+$ 
into $H$-orbits of leaves coincides 
with Lusztig's partition \cite{\Lus} of $G/P_J^+$ 
into locally closed subvarieties which are isomorphic to 
intersections of dual Schubert cells in the full flag variety
$G/B^+$. More precisely, 
$\SS^J_{w_1, w_2}= \P^J_{w_1, w_2', w_2''}$ (in the notation of
\cite{\Rie, Section 5}) where 
$w_2'$ and $w_2''$ are the unique elements in $\WJmin$ and 
$W_J$ such that $w_2 = w_2' w_2''$. This follows, e.g., from the 
discussion in Rietsch \cite{\Rie, Section 5} and \S 1.5 below. 

We also derive explicit formulas for the 
restrictions of the Poisson structures $\pi_J$ to the open 
$B^-$-orbit in Hermitian symmetric spaces of compact type 
for classical groups
and show that in all cases those are the quasiclassical
limits of classes of quadratic algebras that 
attracted a lot of attention in the theory of quantum groups
(see Section 5).
While those classes were previously studied 
case by case, our work conceptually unifies them.
Finally, in the exceptional case
$E_6$ we find a new interesting quadratic Poisson 
structure on a 16 dimensional affine space, 
related to a half-spin representation of $\so_{10}$.

\definition{0.5}
We offer two different proofs of the first part of Theorem 0.4.
The first one, which appears in Section 1, uses the theory of Poisson 
homogeneous spaces. The second one (see Section 3) 
is more geometric and is based upon the notion of Dirac--Poisson 
submanifolds \cite{\CF, \X}. The latter are submanifolds 
of a Poisson manifold $(M, \Pi)$ which in general might not 
be Poisson submanifolds but have the property that their 
symplectic leaves are exactly the connected components
of their intersections with the symplectic leaves of 
$(M, \Pi)$. The second approach uses the idea of 
{\it{weak splittings of surjective Poisson submersions}},
developed in Section 3.
Briefly, if $p : (M, \Pi) \ra (N, \pi)$ is such a 
submersion, a weak splitting of it is a partition 
of $N = \bigsqcup_{\al \in A} N_\al$ into complete Poisson 
submanifolds and liftings $i : N_\al \ra M$ of the 
restrictions of the submersions $p|_{p^{-1}(N_\al)}$
such that  
$(i_\al(N_\al), i_{\al *} (\pi|_{N_\al}))$ 
are Dirac--Poisson submanifolds of $(M, \Pi)$. If
such a weak splitting exists, then the symplectic leaves
of the base $(N, \pi)$ are just the connected components
of the inverse images under the maps $i_\al$ of the 
symplectic leaves of $(M, \Pi)$. It is interesting to 
note that in the present situation the needed partition 
of $(G/P_J^+, \pi_J)$ is exactly the partition
into Schubert cells. The second proof of Theorem 0.4
grew out from our attempt to understand geometrically the observation
\cite{\BGY, Remark 3.10}.
\enddefinition

\definition{0.6}    
As was noted earlier, many interesting quadratic algebras
are quantizations of the algebras of functions on particular 
Schubert cells in particular flag varieties.
The idea of weak splittings of surjective Poisson 
submersions also suggests that the primitive 
ideals of those algebras can be obtained as push forwards 
from the primitive ideals of localizations of quotients 
of the quantized algebras of functions on simple groups
under non-algebra (!) maps. We plan to return to this 
in a forthcoming publication. 

In the case when the parabolic subgroup $P_J^+$ 
has abelian unipotent radical, 
there exists a real form $G_0$ of $G$
for which $G_0 \cap L_J$ is a maximal compact 
subgroup of both $G_0$ and $L_J$. In \cite{\RRS},
 an explicit order-reversing 
bijection is constructed between the $L_J$-orbits on $G/P_J^+$ 
and the $G_0$-orbits on $G/P_J^+$
which were studied in great detail 
in the framework of symmetric spaces by 
Wolf \cite{\Wone}, Takeuchi \cite{\Ta}, and 
others. It is interesting to understand whether 
this bijection can be further refined 
to a bijection between torus orbits of symplectic 
leaves of the Poisson structure $\pi_J$ and
orbits of leaves of a real Poisson structure
on $G/P_J^+$, e.g. the one of Foth and Lu \cite{\FL}.
 
Finally, let us note that parts (i) and (ii) of 
Theorem 0.4 show that each intersection 
of dual Schubert varieties in $G/P_J^+$ 
(i.e., each Richardson variety \cite{\BL} in $G/P_J^+$)
is the Zariski closure of a single $H$-orbit
of symplectic leaves. This suggests the possibility 
to construct explicit degenerations of Richardson 
varieties by deforming algebraically the 
Poisson structure $\pi_J$ and looking at how 
torus orbits of symplectic leaves deform.
This could provide a Poisson geometric approach 
to Schubert calculus.
\enddefinition

\definition{0.7} We conclude the introduction with some notation 
to be used later in the paper.

If $\{\al_1, \ldots, \al_N\}$ denotes the set of positive simple
roots of the reductive Lie algebra $\g$ and $\ga= \sum_i n_i \al_i$ 
is an arbitrary root of $\g$,
then the {\it $\al_j$ height\/} of $\ga$ is defined by
$n_{\al_j}(\ga)= n_j$. The {\it support\/} of $\ga$ is defined by
$$
\supp \ga = \{ \al_j \mid 1 \leq j \leq N, \; \; 
n_{\al_j}(\ga) \ne 0 \}.
$$

If $V$ is a vector space, then for a given $\pi \in \wedge^2 V$ 
we will use the standard notation for the linear map
$$
\pi^\sharp : V^* \ra V, \quad \pi^\sharp(\xi) = \xi \rfloor \pi.
\tag 0.8
$$
For a subspace $V_1 \subset V$, we set
$$
(V_1)^0 = \{ \xi \in V^* \mid \xi(v) = 0 \text{\;for all\;} v \in V_1 \}.
\tag 0.9
$$
In particular, these notations will be used for 
a Poisson bivector $\pi \in (M, \wedge^2 TM)$, in 
which case $\pi^\sharp : T^*M \ra TM$ is a bundle map.
A submanifold $X$ of $(M, \pi)$ will be called a
{\it complete Poisson submanifold\/} if it is 
stable under all Hamiltonian flows, that is, 
if it is union of symplectic leaves. 

Given an algebraic group $G$ and an element $g \in G$, 
we will denote by $\Ad_g$ both the conjugation action 
of $g$ on $G,$ $\Ad_g(h) = g h g^{-1}$, and the 
adjoint action of $g$ on $\Lie G$. As usual, $\ad_x$ 
will be used for the adjoint action of $\Lie G$ on itself.

As in Part I,
we will use the following convention to distinguish
between double cosets and orbits in homogeneous spaces. For two subgroups
$C$ and $D$ of a group $G$ and an element $g \in G:$ (1) The notation
$CgD$ will denote the double coset of $g$ in $G$, and (2) the  notation
$C. gD$ will denote the $C$-orbit of $gD \in G/D$.

Finally, for a locally closed subvariety $Y$ of an algebraic 
variety $X$ and
a subset $Z$ of $Y$, $\Cl_Y(Z)$ will denote the Zariski
closure of $Z$ in $Y$.
\enddefinition


\head 1. Torus orbits of symplectic leaves in flag varieties
via Poisson homogeneous spaces \endhead

\definition{1.1}
In this section we study the geometry of the Poisson structure
$\pi_J$ on the flag variety $G/P_J^+$ by the techniques
of Poisson homogeneous spaces, cf\. \cite{\BGY, Section 1}. 

First we introduce some more notation to be used in the rest 
of the paper and recall basic facts on minimal and maximal 
length representatives for cosets in Weyl groups.
Denote the set of all roots of $(\l_J, \h)$ by 
$\De_J$ and the corresponding subset of positive roots by 
$\De_J^+ = \De_J \cap \De^+$. Recall that the 
(unique) minimal length and maximal length representatives $w$
of a coset in $W/W_J$ are characterized 
respectively by 
$$
w(\De_J^+) \subset \De^+
\tag 1.1
$$
and
$$
w(\De_J^+) \subset -\De^+. 
\tag 1.2
$$
Recall that there are two natural bijection between
$\WJmin$ and $\WJmax$ given by 
$$
v \in \WJmin \mapsto v \wJ \in \WJmax \quad \quad
\text{and} \quad \quad
v \in \WJmin \mapsto w_\ci v \in \WJmax
\tag 1.3
$$ 
where $w_\ci$ and $\wJ$ are the longest elements of 
$W$ and $W_J$.
\enddefinition

\definition{1.2}
Recall that, in the terminology of \cite{\BGY, \S 1.3-1.4},
the standard Poisson algebraic group $(G, \pi_G)$ is a 
part of the algebraic Manin triple 
$(G \times G, G_{\diag}, F)$ where $G_{\diag}$ denotes the 
diagonal subgroup of $G \times G$ and 
$F$ is the dual Poisson algebraic group, given by  
$$
F= \{ (h u^+, h^{-1} u^-) \mid h \in H,\; u^\pm \in U^\pm \}.
\tag 1.4
$$
On the level of Lie algebras, we have the standard Manin triple
$(\g \oplus \g, \g_{\diag}, \g^*)$ where $\g^* = \Lie F$ and
the bilinear form on $\g \oplus \g$ is 
$$ \langle (x_1, x_2), (y_1, y_2) \rangle =
\langle x_1, y_1 \rangle -
\langle x_2, y_2 \rangle,
\tag 1.5
$$
in terms of the nondegenerate bilinear form 
$\langle.,.\rangle$ on $\g$ fixed in \S 0.2.
\enddefinition

\proclaim{1.3. Proposition} {\rm(a)} The orthogonal complement
of $\p_J^+$ in the dual Lie bialgebra $\g^*$
for the standard Lie bialgebra structure on $\g$ 
is $(\p_J^+)^\perp = \n_J^+ \oplus \{0\}$.

{\rm(b)} The standard parabolic subgroup $P_J^+$ of $G$
is a Poisson algebraic subgroup for the standard 
Poisson structure $\pi_G$.

{\rm(c)} The pair $\left( G/P_J^+, \eta_{J*}(\pi_G) \right)$ is
a Poisson homogeneous space for the standard Poisson 
algebraic group $(G, \pi_G)$. The Poisson structure 
$\eta_{J*}(\pi_G)$ is equal to $-\chi(r_\g)$ and
will be denoted by $\pi_J$ for brevity.

{\rm(d)} The Drinfeld 
Lagrangian subalgebra {\rm{(}}cf\. \cite{\Dri}{\rm{)}}
of the base point 
$x_1^J = e P_J^+$ of the Poisson homogeneous space in {\rm(c)} 
is 
$$
\ol{\l}_J =
\{(l+n_1^+, l+n_2^+) \mid l \in \l_J,\; n^+_i \in \n_J^+ \}
\subset \g \oplus \g \cong D(\g).
\tag 1.6
$$
It is the tangent Lie algebra of the connected algebraic subgroup
$$
\ol{L}_J =
\{(l u_1^+, l u_2^+) \mid l \in L_J,\; u^+_i \in U_J^+ \} 
\tag 1.7
$$
of $G \times G$; in particular,
$\left( G/P_J^+, \pi_J \right)$
is an algebraic Poisson homogeneous space for the standard 
Poisson algebraic group $(G, \pi_G)$ in the terminology 
of \cite{\BGY, Definition 1.7}.
\endproclaim

The proof of Proposition 1.3 is analogous to that of 
\cite{\BGY, Proposition 3.2} and will be omitted. 

Below we will need the following well known Lemma.

\proclaim{1.4. Lemma} All $B^-$-orbits on the flag variety 
$G/P_J^+$ are parametrized by $\WJmax$ by $w \mt B^- x_w^J.$ 
Moreover $B^- x_w^J$ is biregularly isomorphic to $U_w^-$, 
recall {\rm(0.7)}, by
$$
u \in U_w^- \mt u x_w^J.
$$  
\endproclaim

For completeness, we will 
sketch the proof of the second part of Lemma 1.4. 
For any $w \in W$, one has $B^-= U^-_w (B^- \cap \Ad_w(B^+))$. 
If $w \in \WJmax$ then (1.2)
implies that $\n^- \cap \Ad_w(\n^-) \subset \Ad_w(\n^-_J)$.
Since $U^-_w$ and $\Ad_w(U_J^-)$ are connected
subgroups of $G$, it follows that $U^-_w \subset \Ad_w(U_J^-)$,
which easily implies the second statement in 
Lemma 1.4.
  
Set 
$$
\Om^J =\{(w_1, w_2) \in \WJmax \times W \mid
w_1 \leq w_2 \}.
\tag 1.8
$$ 

\proclaim{1.5. Theorem} There are finitely many $H$-orbits 
of symplectic leaves on $(G/P_J^+, \pi_J)$, bijectively parametrized by 
$\Omega^J$, and 
all of them 
are smooth irreducible locally closed subvarieties of $G/P_J^+$.
The $H$-orbit of leaves corresponding to $(w_1, w_2) \in \Om^J$
is explicitly given by
$$
\SS_{w_1,w_2}^J = \U_{\dw_1, w_2} . P_J^+,
\tag 1.9
$$  
recall {\rm(0.6)}, and is biregularly isomorphic to 
the intersection $\B_{w_1,w_2}= B^- x_{w_1} \cap B^+ x_{w_2}$
of dual Schubert
cells in the full flag variety $G/B^+$.

In particular, the $H$-orbits of symplectic
leaves on $(G/B^+, \pi)$ are exactly the 
intersections of dual Schubert cells
$\B_{w_1,w_2}$, indexed by pairs 
$(w_1, w_2) \in W \times W$ 
such that $w_1 \leq w_2$. 
\endproclaim

It is easy to see that, while $\U_{\dw_1, w_2}$ 
depends on the representative of $w_1$ in the normalizer of 
$H$, the set $\SS_{w_1, w_2}^J$ does not. 

\demo{Proof} It follows from \cite{\BGY, Theorem 1.10} 
that the $H$-orbits
of symplectic leaves in $(G/P_J^+, \pi_J)$
are smooth locally closed subsets of $G/P_J^+$. 
Moreover the same Theorem implies that they are 
exactly the irreducible components 
of the inverse images under the map
$$
\De : G/P_J^+ \hra (G\times G)/\ol{L}_J, \quad 
\De(g P_J^+)= (g,g) \ol{L}_J
\tag 1.10
$$
of the $(B^+ \times B^-)$-orbits on $(G\times G)/\ol{L}_J$.
The map $\De$ is an embedding because
$$
G_{\diag} \cap \ol{L}_J = (P_J^+)_{\diag}.
$$ 
(Recall that $G_{\diag}$ denotes the diagonal subgroup
of $G \times G$.)

Applying \cite{\BGY, Theorem 8.1}, 
we see that each such orbit passes
through the coset of $(\dw_2, \dw_1)$ for some  
$(w_1, w_2) \in \WJmin \times W$ and all such orbits 
are distinct. Because 
$(\dw^J_\ci, \dw^J_\ci) \in \ol{L}_J$, using (1.3) 
we obtain that 
$$
\{(\dw_2, \dw_1) \ol{L}_J \mid w_1 \in \WJmax,\; w_2 \in W \}
$$
is a complete, 
irredundant set of representatives
for all $(B^+ \times B^-)$-orbits on $(G\times G)/\ol{L}_J$.
Note that 
$$
\De^{-1} \left( (B^+ \times B^-) (\dw_2, \dw_1) \ol{L}_J \right)
\subseteq B^- x_{w_1}^J.
$$
Thus, 
$\De^{-1} \left( (B^+ \times B^-) (\dw_2, \dw_1) \ol{L}_J \right)$
consists of all points $u^- x_{w_1}^J$, with $u^- \in U^-_w$, for which
there  exist $b^\pm \in B^\pm$, $u^+_1,u^+_2 \in U_J^+$, and $l
\in L_J$ such that   
$$
u^- \dw_1 = b^- \dw_1 l u^+_2 = b^+ \dw_2 l u^+_1 .
\tag 1.11
$$
From the first equality we get that 
$l u^+_2 = \Ad^{-1}_{\dw_1}((b^-)^{-1}u^-) \in \Ad^{-1}_{w_1}(B^-)$.
Since $w_1 \in \WJmax$ we have $\Ad_{\dw_1}(u^+_2) \in B^-$, thus
$l \in L_J \cap \Ad_{w_1}(B^-)$. For arbitrary
$l \in L_J \cap \Ad_{w_1}(B^-)$
and $u^- \in U^-_w$ there exist $b^- \in B^-$ and 
$u_2^+ \in U_J^+$ that satisfy the first equality in (1.11). 
Then the second equality in (1.11) implies  
$$
\De^{-1} \left( (B^+ \times B^-) (\dw_2, \dw_1) \ol{L}_J \right)
= \left( U^-_{w_1} \dw_1 \cap B^+ w_2 
(L_J \cap \Ad_{w_1}^{-1}(B^-)) U_J^+
\right) . P_J^+.  
$$
From (1.2) one has 
$\Ad_{w_1}^{-1}(B^-)=(L_J \cap B^+)(U_J^- \cap\Ad^{-1}_{w_1}(B^-)).$
Thus $L_J \cap \Ad^{-1}_{w_1}(B^-)= L_J \cap B^+$ and
$$
\De^{-1} \left( (B^+ \times B^-) (\dw_2, \dw_1) \ol{L}_J \right)
= \left( U^-_{w_1} \dw_1 \cap B^+ w_2 B^+ \right) . P_J^+
= \U_{\dw_1, w_2}. P_J^+.  
$$
Due to Deodhar's theorem \cite{\Deo, Corollary 1.2}, 
$\U_{\dw_1, w_2}$ is non-empty if and only if $w_2 \geq w_1$.
On the other hand if $w_2 \geq w_1$ then $\SS^J_{w_1, w_2}$
is irreducible and biregularly isomorphic to $\U_{\dw_1, w_2}$
and $\B_{w_1,w_2}$ because of Lemma 1.4 and 
\cite{\BGY, Theorem 2.4}.
\qed\enddemo

\definition{1.6} Next we describe the relation between 
the Poisson structures on torus orbits of leaves in different 
flag varieties for a fixed reductive group $G$. For a subset 
$J \subseteq \Ga$, denote by
$$
\mu_J : G/B^+ \ra G/P_J^+
\tag 1.12
$$
the natural projection. The composition 
$G @>{\eta}>> G/B^+ @>{\mu_J}>> G/P_J^+$
coincides with $\eta_J$,
recall (0.3). Since the surjective maps $\eta$ and $\eta_J$
are Poisson, the projection $\mu_J$ is Poisson as well.

For all $J \subseteq \Ga$ and $w \in W$, the sets
$B^\pm x_w^J$ are 
complete Poisson subvarieties of 
$(G/P_J^+, \pi_J)$. This follows from the facts 
that $(G/P_J^+, \pi_J)$ is a quotient of
$(G, \pi_G)$ and $B^\pm$ are Poisson algebraic
subgroups of $(G, \pi_G)$. (The statement is
also a simple corollary of Theorem 1.5.)
\enddefinition

\proclaim{Proposition} For every subset $J \subseteq \Ga$
and every $w \in \WJmax$, the projection $\mu_J$ restricts 
to an isomorphism of Poisson varieties
$$
\mu_J|_{B^- x_w} : (B^- x_w, \pi) @>{\cong}>> 
(B^- x_w^J, \pi_J).
$$
In particular, for $(w_1, w_2) \in \Om^J$ this restricts to the 
Poisson isomorphism
$$
\mu_J|_{\SS_{w_1,w_2}} : (\SS_{w_1,w_2}, \pi) @>{\cong}>>
(\SS^J_{w_1, w_2}, \pi_J)
$$
where for simplicity we set 
$\SS_{w_1,w_2}:= \SS^\emptyset_{w_1,w_2}$.
\endproclaim

\demo{Proof} Lemma 1.4 guarantees that 
$\mu_J$ is an isomorphism between the affine spaces
$B^- x_w \subset G/B^+$ and $B^- x_w^J \subset G/P^+_J$. 
It was shown above that 
$B^- x_w$ and $B^- x_w^J$ are complete Poisson 
subvarieties of $(G/B^+, \pi)$ and 
$(G/P_J^+, \pi_J)$, respectively.
The first statement 
now follows from the fact that $\mu_J$ is 
a Poisson mapping and the second one 
is a corollary of the first. 
\qed
\enddemo

\proclaim{1.7. Remark} One can establish first 
Proposition {\rm 1.6} and then deduce the general case 
of Theorem {\rm 1.5} from the case $J = \emptyset.$
The proof of Theorem {\rm 1.5} in this special case 
is easier than the general case.
This gives a somewhat simpler proof of Theorem {\rm 1.5}. 
\endproclaim

Finally we describe the Zariski closures of the torus orbits 
of symplectic leaves $\SS^J_{w_1, w_2}$
of $(G/P_J^+, \pi_J)$.

\proclaim{1.8. Theorem} For all $(w_1, w_2) \in \Om^J$, 
the Zariski closure of the $H$-orbit
of symplectic leaves $\SS^J_{w_1, w_2}$
in $(G/P_J^+, \pi_J)$ is equal to the 
union of $\SS^J_{v_1, v_2}$ over those 
$(v_1, v_2) \in \Om^J$ for which there exists
$z \in W_J$ such that 
$w_1 \leq v_1 z$ and 
$w_2 \geq v_2 z$.    

\endproclaim 

Note that in the special case of the generalized full
flag variety $G/B^+$, Theorem 1.8 reduces to the well 
known fact that for the Richardson varieties,
$$ 
\ol{B^- x_{w_1} \cap B^+ x_{w_2}} =
\ol{B^- x_{w_1}} \cap \ol{B^+ x_{w_2}}.
$$

The order relation between the nonnegative parts of 
$\SS^J_{w_1, w_2}$ in the sense of Lusztig \cite{\Lus}
was also (independently) studied by Rietsch \cite{\Rie}
and Williams \cite{\Wi}.

\demo{Proof} By modifying 
Springer's Lemma 2.2 \cite{\Sp}
or directly using \cite{\LY, Corollary 3.6},
one gets that for all 
$(w_1, w_2) \in \WJmin \times W$ inside 
$(G \times G)/\ol{L}_J$: 
$$
\aligned
\ol{(B^+ \times B^+).(w_2, w_1) \ol{L}_J}=
\bigsqcup \{ & (B^+ \times B^+).(v_2, v_1) \ol{L}_J \mid
(v_1, v_2) \in \WJmin \times W,
\\ &\exists z \in W_J \; \; 
\text{such that} \; \; 
w_1 \geq v_1 z,\; w_2 \geq v_2 z \}.
\endaligned
\tag 1.13 
$$
Acting on (1.13) by $(e,\dw_\ci)$, and using that 
$\Ad_{w_\ci}(B^+) = B^-$ and that for 
$w, v \in W$, we have $w \leq v$ if and only if
$w_\ci w \geq w_\ci v$, we see that for 
$(w_1, w_2) \in \WJmin \times W$,
$$
\aligned
\ol{(B^+ \times B^-).(w_2, w_\ci w_1) \ol{L}_J}=
\bigsqcup \{ & (B^+ \times B^-).(v_2, w_\ci v_1) \ol{L}_J \mid
(v_1, v_2) \in \WJmin \times W, \\
&\exists z \in W_J \; \; 
\text{such that} \; \;
w_\ci w_1 \leq w_\ci v_1 z,\; w_2 \geq v_2 z\}.
\endaligned
$$
We deduce from the second correspondence 
between $\WJmin$ and $\WJmax$ in (1.3)
that for all $(w_1, w_2) \in \WJmax \times W$,
$$
\aligned
\ol{(B^+ \times B^-).(w_2, w_1) \ol{L}_J}=
\bigsqcup \{ & (B^+ \times B^-).(v_2, v_1) \ol{L}_J \mid
(v_1,v_2) \in \WJmax \times W, \\
&\exists z \in W_J \; \; 
\text{such that} \; \;
w_1 \leq v_1 z,\; w_2 \geq v_2 z\}. 
\endaligned
\tag 1.14
$$

Now we apply \cite{\BGY, Lemma 2.6} for $Y$ equal to the 
image of the embedding $\De$ from (1.1) 
and the stratification of $(G\times G)/\ol{L}_J$
by $(B^+ \times B^-)$-orbits. Note that $\De(G/P_J^+)$
intersects each $(B^+ \times B^-)$-orbit transversally since the 
diagonal of $\g \oplus \g$ and $\b_+ \oplus \b_-$
span $\g \oplus \g$. It follows that
$$
\multline
\Cl_{\De(G/P_J^+)} 
\left( \De(G/P_J^+) \cap 
((B^+ \times B^-).(w_2, w_1) \ol{L}_J) \right) 
\\
=\De(G/P_J^+) \cap 
\ol{(B^+ \times B^-).(w_2, w_1) \ol{L}_J}.
\endmultline
\tag 1.15 
$$
Recall from the proof of Theorem 1.5 that for 
$(w_1, w_2) \in \WJmax \times W$,
$$
\De(G/P_J^+) ((B^+ \times B^-).(w_2, w_1) \ol{L}_J)
\ne \emptyset \quad \text{if and only if} \quad w_1 \leq w_2,
$$
in which case
$$
\SS^J_{w_1, w_2} = \De^{-1}\left( 
(B^+ \times B^-).(w_2, w_1) \ol{L}_J \right).
$$
The Lemma now follows from (1.14), (1.15), and the 
fact that $\De$ is an embedding.
\qed
\enddemo
The Zariski closures of the $H$-orbits of symplectic leaves 
inside each Schubert cell in $(G/P_J^+, \pi_J)$ have a 
particularly simple form.

\proclaim{1.9. Proposition} For each $w \in \WJmax$
the following hold.

{\rm(i)} The $H$-orbits of symplectic leaves in 
the Schubert cell 
$(B^- x_w^J \subset G/P_J^+, \pi_J|_{B^- x_w^J})$
are parametrized by 
$W^{\geq w} = \{ w_2 \in W \mid w_2 \geq w \}$ by 
$$
w_2 \in W^{\geq w} \mapsto 
\SS^J_{w, w_2} = \U_{w, w_2} . P_J^+.
$$

{\rm(ii)} For each $w_2 \in W^{\geq w}$,
the Zariski closure of $\SS^J_{w, w_2}$ inside the Schubert
cell $B^- x_w^J$ consists of all 
$\SS^J_{w, v_2}$ for $v_2 \in W$ such that 
$w\leq v_2 \leq w_2$.
\endproclaim 

This proposition generalizes Theorems 3.9 and 3.13
from \cite{\BGY}.

\demo{Proof} The first part follows from 
$$
\SS^J_{w, w_2} = \U_{w, w_2} . P_J^+ \subset (U_w^- w). P_J^+=
B^- x_w^J 
$$
and the fact that the Schubert cells $B^- x_w^J$ for
$w \in \WJmax$ partition $G/P_J^+$.

It is easy to deduce the second part from Theorem 1.8
but we offer a direct proof which better explains 
the result.

In terms of the isomorphism
$\mu_J|_{B^- x_w} : B^- x_w \ra B^- x_w^J$ (cf\. \S 1.6),
$\SS^J_{w, w_2}$ is given by
$$
\SS^J_{w, w_2} = \left( \mu_J|_{B^- x_w} \right)^{-1}
(B^- x_w \cap B^+ x_{w_2}). 
$$
Applying \cite{\BGY, Lemma 2.6} for the partition 
$X= G/B^+ = \bigsqcup \Sb w_1 \in W \endSb B^+ x_{w_1}$ and
$Y = B^- x_w$, and using the standard formulas
for closures of Schubert cells, leads to
$$ 
\Cl_{B^- x_w} (B^- x_w \cap B^+ x_{w_2}) = 
\bigsqcup \Sb v_2 \in W, w \leq v_2 \leq w_2 \endSb 
(B^- x_w \cap B^+ x_{v_2}). 
$$
The second part of the Proposition follows from 
this, applying once again the 
isomorphism $\mu_J|_{B^- x_w}$.
\qed\enddemo


\head 2. Weak splittings of surjective Poisson submersions 
\endhead

\definition{2.1} First we recall the definition of a
Poisson--Dirac submanifold of a Poisson manifold, 
given by Crainic and Fernandes in \cite{\CF, Section 9}.
\enddefinition

\proclaim{Definition} Assume that $(M, \Pi)$ is a 
smooth {\rm(}real or complex{\rm)} Poisson manifold.
A submanifold $X$ of $M$ is called a 
\underbar{Poisson--Dirac submanifold} if the following two
conditions are satisfied:

{\rm(i)} For each symplectic leaf $S$ of $(M, \Pi)$, 
the intersection $S\cap X$ is clean {\rm(}i.e., it is smooth and 
$T_x (S \cap X)= T_x S \cap T_x X$
for all $x \in S \cap X${\rm)} and
$S\cap X$ is a symplectic submanifold of 
$(S, (\Pi|_S)^{-1})$.

{\rm(ii)} The family of symplectic structures 
$(\Pi|_S)^{-1}|_{S\cap X}$ is induced by a smooth Poisson 
structure $\pi$ on $X$.
\endproclaim
Here and below, for a nondegenerate Poisson structure $\pi_0$ 
we denote by $(\pi_0)^{-1}$ the corresponding symplectic form.

In the setting of the above Definition, 
the symplectic leaves of $(X, \pi)$ are 
exactly the connected components of the 
intersections of the symplectic leaves 
of $(M, \Pi)$ with $X$. 

\definition{2.2} The following simple criterion was 
proved in \cite{\CF}.
\enddefinition

\proclaim{Proposition} Assume that $(M, \Pi)$ is 
a Poisson manifold and that $X$ is a submanifold
for which there exists a subbundle $E$ of $T_X M$
such that 

{\rm(i)} $T_X M = T X \oplus E$ and

{\rm(ii)} the restriction of the Poisson tensor $\Pi$ to $X$ 
splits as
$$
\Pi|_X = \pi + \pi_E
$$
for some smooth bivector fields $\pi \in \Ga(X, \wedge^2 TX)$
and $\pi_E \in \Ga(X, \wedge^2 E)$.

Then $X$ is a Poisson--Dirac submanifold of $(M, \Pi)$
and the induced Poisson structure on it coincides with $\pi$.
\endproclaim

Crainic and Fernandes call submanifolds satisfying the 
conditions of Proposition 2.2 {\it{Poisson--Dirac submanifolds
admitting a Dirac projection}}. Earlier, Xu \cite{\X}
investigated such submanifolds 
with an extra property, namely that $E^0$ 
is a Lie subalgebroid of $T^*M,$ equipped with 
the standard cotangent bundle algebroid structure,
recall (0.9).

\proclaim{2.3. Definition} Assume that $(M, \Pi)$
and $(N, \pi)$ are Poisson manifolds and that 
$p : (M, \Pi) \ra (N, \pi)$ is a surjective Poisson 
submersion. A \underbar{weak splitting} of $p$ is a partition
$$
N = \bigsqcup_{\al \in A} N_\al
\tag 2.1
$$
of $(N, \pi)$ into complete Poisson submanifolds
such that for each $\al \in A$, there exists
a smooth lifting $i_\al : N_\al \ra M$
{\rm(}of $p|_{p^{-1}(N_\al)} : p^{-1}(N_\al) \ra N_\al${\rm)}
with the properties:

{\rm(i)} $i_\al(N_\al)$ is a Poisson--Dirac
submanifold of $(M, \Pi)$ and

{\rm(ii)} the induced Poisson structure on $i_\al(N_\al)$
is $i_{\al *}(\pi|_{N_\al})$. 
\endproclaim

\noindent
Note that $i_\al$ {\it{is not required to be 
a Poisson map}}.
An important special case is illustrated in algebraic terms
in Proposition 2.6.

\proclaim{2.4. Remark} If a surjective Poisson 
submersion $p : (M, \Pi) \ra (N, \pi)$ admits a 
weak splitting as in Definition {\rm2.3}, then the 
symplectic foliation of $(N, \pi)$ is easily
described in terms of the symplectic foliation
of $(M, \Pi)$. Namely, each symplectic leaf of 
$(N, \pi)$ lies entirely in one of the submanifolds
$N_\al$ and is of the type
$i_\al^{-1}(S \cap i_\al(N_\al))^\ci$
where $S$ is a symplectic leaf of $M$ and $(S \cap i_\al(N_\al))^\ci$ is
a  connected component of
$(S \cap i_\al(N_\al))$.
\endproclaim

\definition{2.5} The following Proposition provides a 
sufficient condition for the condition 
(ii) in Definition 2.3 which is easier to check.  
\enddefinition

\proclaim{Proposition} Assume that 
$p : (M, \Pi) \ra (N, \pi)$ is a surjective Poisson 
submersion. Let 
$$
N = \bigsqcup_{\al \in A} N_\al
$$
be a partition of $(N, \pi)$ into complete Poisson 
submanifolds such that for each $\al \in A$, 
there exists a smooth lifting $i_\al : N_\al \ra M$
{\rm(}of $p|_{p^{-1}(N_\al)} : p^{-1}(N_\al) \ra N_\al${\rm)}
whose image is a Poisson--Dirac submanifold admitting
a Dirac projection with respect to a subbundle
$E_\al$ of $T_{i_\al(N_\al)}M$,
cf\. Proposition {\rm2.2}. If the
$E_\al$ are tangent to the fibers of $p$
{\rm(}i.e., $E_\al$ contain the tangent spaces to the 
fibers of $p${\rm)}, then 
the condition {\rm(ii)} in Definition {\rm2.3} is satisfied
and the families $\{N_\al\}_{\al \in A}$,
$\{i_\al\}_{\al \in A}$ provide a weak splitting of 
$p$.
\endproclaim

\demo{Proof} For $m \in i_\al(N_\al)$ denote
the fiber of $p$ through $m$ by $F_m = p^{-1}(p(m))$.
Since $p : (M, \Pi) \ra (N, \pi)$ is Poisson and
$i_\al$ is a lifting of
$p|_{p^{-1}(N_\al)} : p^{-1}(N_\al) \ra N_\al$
we have
$$
T_m M = T_m (i_\al(N_\al)) \oplus T_m F_m
\tag 2.2
$$
and 
$$
\Pi_m- i_{\al*}(\pi|_{N_\al})_m \in T_m F_m \wedge T_m M.
\tag 2.3
$$
On the other hand, 
the fact that $i_\al(N_\al) \subset M$ satisfies
the conditions of Proposition 2.2 implies
$$
\aligned
\Pi|_{i_\al(N_\al)} &=  \pi_\al + \pi_{E_\al}, \\
 &\; \text{for some}
\; \pi_\al \in 
\Ga(i_\al(N_\al), \; \wedge^2 T (i_\al(N_\al))), \; 
\pi_{E_\al} \in 
\Ga(i_\al(N_\al), \wedge^2 E_\al).
\endaligned \tag 2.4
$$ 
Putting together (2.2)--(2.4) and $T_m F_m \subseteq E_\al(m)$
gives $\pi_\al = i_{\al*}(\pi|_{N_\al})$ 
which is exactly the condition (ii) in Definition 2.4
(taking into account Proposition 2.2).
\qed\enddemo

\definition{2.6} If $p : (M, \Pi) \rightarrow (N, \pi)$ is a Poisson
map,
then the pull back $p^* : (C^\infty(N), \{.,.\}_\pi) \rightarrow
(C^\infty(M), \{.,.\}_\Pi)$ is a homomorphism of Poisson algebras
and turns $C^\infty(M)$ into a module for the Poisson algebra
$(C^\infty(N), \{.,.\}_\pi)$. The following Proposition provides an
algebraic characterization of an important special case of weak
splittings of surjective Poisson submersions. Its proof is simple and
will be left to the reader.
\enddefinition
 
\proclaim{Proposition} Assume that $p : (M, \Pi) \rightarrow (N, \pi)$
is a surjective Poisson submersion and $i : N \rightarrow M$
 a smooth lifting of $p$. Denote by $Tp$ the subbundle of $TM$
whose fibers are the tangent spaces to the fibers of $p$. Then the
trivial partition of $N$ with one stratum and the map
$i : N \rightarrow M$ provide a weak splitting of $p$ such that
$\Pi|_{i(N)} \in \wedge^2 T i(N) \oplus \wedge^2 T p|_{i(N)}$
if and only if  $i^* : (C^\infty(M), \{.,.\}_\Pi) \rightarrow
(C^\infty(N), \{.,.\}_\pi)$ is a morphism of $(C^\infty(N),
\{.,.\}_\pi)$
modules. 
\endproclaim

\definition{2.7} All weak splittings of surjective Poisson 
submersions considered in this paper will be in the category
of (complex) quasiprojective Poisson varieties. This means that in 
the setting of Definition 2.1, we require $X$ to be a (smooth) 
locally closed subvariety of the smooth quasiprojective Poisson 
variety $M.$ In Proposition 2.2, we require $E$ to be an algebraic 
subbundle of $T_X M$. Finally, in the algebraic setting,
in Definition 2.3 we require (2.1) to be an algebrogeometric
stratification of $M$ (in the sense of \cite{\BGY, \S 0.8}) 
and the maps $i_\al$ to be algebraic.
\enddefinition


\head 3. Weak splittings of surjective Poisson 
submersions for flag varieties 
\endhead
\definition{3.1} Since the Poisson structure 
$\pi_G$ vanishes on the maximal torus $H$ of $G$,
the left and right regular 
actions of $H$ on $G$ preserve it, and thus
$\pi_G$ descends to a Poisson structure 
on $G/H$. One can see this in another way: 
$H$ is a Poisson 
algebraic subgroup of $(G, \pi_G)$, thus
$\pi_G$ descends to a Poisson structure on the 
homogeneous space $G/H$ because of 
\cite{\BGY, Theorem 1.8}. Denote the 
standard projection by
$$
\tau : G \ra G/H, \quad \text{and set} \quad
\pi_{G/H}= \tau_*(\pi_G).
\tag 3.1
$$
It is clear that the projections
$$
\nu_J : (G/H, \pi_{G/H}) \ra (G/P_J^+, \pi_J), 
\quad \nu_J(g H) = g P_J^+
\tag 3.2
$$ 
are surjective Poisson submersions. 
For brevity, set $\nu=\nu_\emptyset$.
Finally, for $w \in W$, set
$$
y_w= w H \in G/H.
\tag 3.3    
$$
The following Theorem contains the main 
result in this Section. Based on it,
a second proof of Theorem 1.5 is given in 
\S 3.9.
\enddefinition

\proclaim{3.2. Theorem} Assume that 
$G$ is an arbitrary complex reductive algebraic group 
and $J$ is a subset of the set of positive simple roots $\Ga$.
The partition into Schubert cells
$$
G/P_J^+ = \bigsqcup \Sb w \in \WJmax \endSb B^- x_w
$$ 
and the morphisms 
$$
i_w^J : B^- x_w^J \ra G/H,
\quad
\text{given by} \quad
i_w^J(u^- x_w^J) \ra u^- y_w, 
\; \text{for} \; u^- \in U^-_{w},
\tag 3.4
$$
provide a weak splitting, in the sense of Definition {\rm2.3}, 
 of the surjective Poisson
submersion $\nu_J : (G/H, \pi_{G/H}) \ra (G/P_J^+, \pi_J)$
{\rm (}recall {\rm(0.5)}, {\rm(0.7)}, and 
{\rm(3.1)--(3.3))}.
\endproclaim

The proof of Theorem 3.2 will be split into several Lemmas.

\definition{3.3} We will use the identification of vector 
spaces
$$
T_e G \oplus T_e^* G = \g \oplus \g^* \cong 
D(\g) \cong \g \oplus \g = T_e G \oplus T_e G.
\tag 3.5
$$
coming from the embeddings of $\g$ and $\g^*$ in the 
double $D(\g) \cong \g \oplus \g$ 
of the Lie bialgebra $\g$, cf\. \S 1.2.
This induces the identifications
$$ 
T_g G \oplus T_g^* G \cong R_g(\g \oplus \g^*) \cong 
R_g(\g \oplus \g) = T_g G \oplus T_g G,
\quad \text{for} \; g \in G.
\tag 3.6
$$
Here and below, $L_g$ and $R_g$ refer to 
the left and right translations $a \mapsto g a$ and
$a \mapsto a g$, for $a \in G$. For simplicity of the 
notation, we denote in the same way the induced
tangent maps $T_a \ra T_{g a}$ and $T_a \ra T_{a g}$. 
In the identifications (3.6),
the tangent and cotangent spaces at $g$ 
correspond respectively to: 
$$
\aligned
T_g G & \cong R_g(\g_{\diag}),
\\
T_g^* G & \cong R_g(\Lie F) =
R_g \{ (h + n^+, - h+ n^-) \mid h \in \h,\; n^\pm \in \n^+ \},
\endaligned
\tag 3.7
$$
recall the notation in \S 1.2. The pairing between them is 
given by (1.5).

The main reason for using the identifications (3.6) 
is that the graph of
$\pi_{G ,g}^\sharp : T_g^* G \ra T_g G$ (cf\. (0.8))
under those identifications corresponds to the Drinfeld
Lagrangian subalgebra \cite{\Dri} $\l_g$ 
of $D(\g)$ for the base point $g$ of $(G, \pi_G)$,
considered as a homogeneous space over itself. 
Since $\pi_G$ vanishes at $e$, one has $\l_e= \Lie F$. Moreover,
because the map $g \mapsto \l_g \subset D(\g)$ is 
$G$-equivariant with respect to the adjoint action
of $G$ on $D(\g) \cong \g \oplus \g$ (see \cite{\Dri}),
one has
$$
\l_g = \Ad_g(\Lie F) = 
\Ad_g \left( \{ (h + n^+, - h+ n^-) \mid h \in \h,\; 
n^\pm \in \n^+ \}\right).
$$
Taking into account $R_g \circ \Ad_g = L_g$, 
this leads us to the following result.
\enddefinition

\proclaim{3.4. Lemma}
In the identification {\rm(3.6)}, the graph of 
$\pi_{G ,g}^\sharp : T_g^* \ra T_g G$, cf\. {\rm(0.8)},
corresponds to the subspace
$$
R_g (\Ad_g F) = L_g (\Lie F) = 
L_g \left( \{ (h + n^+, - h+ n^-) \mid h \in \h,\; 
n^\pm \in \n^+ \}\right)
\subset T_g G \oplus T_g G.
$$
\endproclaim
We will further need the following well known result.

\proclaim{3.5. Lemma} If $V_1$ and $V_2$ are subspaces 
of a finite dimensional vector space $V$ and 
$\pi \in \wedge^2 V$, then 
$$
\pi^\sharp(V_1^0) \subseteq V_2 \quad \iff \quad 
\pi \in \wedge^2 V_1 + \wedge^2 V_2.
$$
\endproclaim

\demo{Sketch of the proof} The standard nondegenerate pairing 
between $V \oplus V$ and $V^* \oplus V^*$ restricts to a 
nondegenerate pairing between $\wedge^2 V$ and $\wedge^2 V^*$.
Then $\pi^\sharp(V_1^0) \subseteq V_2$ if and only if 
$\pi \in (V_1^0 \wedge V_2)^0$. One easily checks
that for all subspaces $V_1$ and $V_2$ of $V$, one has
$(V_1^0 \wedge V_2)^0 = \wedge^2 V_1 + \wedge^2 V_2$.
\qed
\enddemo

\definition{3.6}
For all $w \in W$, define the following algebraic subbundles
$\wt{E}_w$ and $E_w$ of $T_{H U^-_w w} G$ and 
$T_{U^-_w y_w} G/H$, respectively:
$$
\aligned
\wt{E}_w(bw) &= L_{bw} (\b^+) + 
R_{bw} (\n^+ \cap \Ad_w(\n^-)) \subset T_{b w} (G), 
\quad \text{for} \; b \in H U^-_w,
\\
E_w & = \tau_* \left( \wt{E}_{w} \right) 
\subset T_{U^-_w y_w} G/H, 
\endaligned
\tag 3.8
$$ 
recall (3.4). It is easy to see that the push-forward in (3.8)
does not depend on the choice of preimage. 
For a subvariety $N \subset G$ such that $N = \tau^{-1} \tau (N)$, 
denote by $T (\tau|_N)$ the bundle over $N$ whose fibers are the 
tangent spaces to the fibers of $\tau|_{N}$:
$$
T( \tau|_{N} ) (g) = L_g(\h) \subset T_g G, 
\quad g \in N.
\tag 3.9
$$ 
The fact that the $E_w$ are algebraic bundles follows 
from $\wt{E}_w \supset T (\tau|_{H U^-_w w})$: 
for all $w \in W$ and $b \in H U^-_w$, one has
$\wt{E}_w (b w) \supset L_{bw}(\b^+) \supset L_{bw} (\h) 
= T( \tau|_{H U^-_w w})(bw)$.
\enddefinition

\proclaim{Proposition} For every $w \in W$, 
the following hold:
\roster
\item 
$T (H U_w^- w) + \wt{E}_w = T_{H U_w^- w} G$,
\item""
$T (H U_w^- w) \cap \wt{E}_w = T(\tau|_{H U_w^- w})$,
\item""
$\pi_{G}|_{H U_w^- w} \in \wedge^2 T (H U_w^- w)
       + \wedge^2 \wt{E}_w$, and
\item
 $T (U_w^- y_w) \oplus E_w = T_{U_w^- y_w} G$,
\item""
$\pi_{G/H}|_{U_w^- y_w} \in \wedge^2 T (U_w^- y_w)
       \oplus  \wedge^2 E_w,$ 
\endroster
recall {\rm(0.4)}, {\rm(3.3)}, {\rm(3.4)}, {\rm(3.8)}, and {\rm(3.9)}. 
\endproclaim

\demo{Proof} Part (2) follows from part (1). 
The following inclusions imply the first statement 
in part (1):
$$
\aligned
T_{bw} (H U^-_w w) + \wt{E}_w(uw) 
 &\supset  R_{bw} \bigl( \b^- \cap \Ad_w(\b^-) + \n^+ \cap \Ad_w(\n^-)
\bigr) + L_{bw} (\b^+)  \\
 &= R_{bw} (\Ad_w(\b^-)) + L_{uw} (\b^+)
= R_{bw} (\Ad_{bw}(\b^-)) + L_{bw} (\b^+)  \\
 &= L_{bw} (\b^- + \b^+) = T_{bw}(G),
\endaligned
$$
where $w \in W$ and $b \in U^-_w$. We used that 
$b \in H U^-_w \subset \Ad_w(B^-)$ and thus 
$\Ad_{bw} (\b^-) = \Ad_{w}(\b^-)$.

It is clear that  
$$
T (H U_w^- w) \cap \wt{E}_w \supset T(\tau|_{H U_w^- w})
\tag 3.10
$$
because $H U_w^- w = U_w^- w H$. Next we show the 
opposite inclusion.
Fix $b \in H U^-_w$.
If $x\in \g$ and $R_{bw}(x) \in T_{bw} (H U_w^- w) \cap \wt{E}_w$, then 
$x = x_1 =\Ad_{bw}(y) + x_2$ for 
some $y \in \b^+$ and $x_1 \in \b^- \cap \Ad_w(\b^-)$, 
and $x_2 \in \n^+ \cap \Ad_w(\n^-)$. So 
$$
\Ad_{bw} (y) = x_1 - x_2 \in 
\b^- \cap \Ad_w(\b^-) + \n^+ \cap \Ad_w(\n^-)
=\Ad_{w}( \b^-) = \Ad_{bw}(\b^-).
$$
Taking into account $\Ad_{bw}(y) \in \Ad_{bw} (\b^+)$,
one obtains $y\in h$ and $x_2 = 0$, and consequently
$x \in \h$. 
This proves the opposite inclusion to (3.10) and 
completes the proof of the second statement in (1).

In the rest of this proof we show the third property in 
part (1). Under the identification (3.6), the tangent space
$T_{bw} (HU^-_w w)$ corresponds to 
$R_{bw}\left((\b^- \cap \Ad_w(\b^-))_{\diag}\right)$
for all $b \in H U^-_w$.
Recall that the images of the tangent and the cotangent 
spaces to $G$ in (3.6) are given by (3.7) and the pairing 
between them is given by (1.5), see \S 3.3. By a direct 
computation, one checks that under the identification (3.6), 
$(T_{bw}(H U^-_w w))^0$ corresponds to 
$R_{bw}((\n^+ \cap \Ad_w(\n^-)) \oplus \n^-)$, see (0.8)--(0.9).
(Here and below, for two subalgebras $\f_1, \f_2 \subseteq \g$, by 
$\f_1 \oplus \f_2$ we denote the canonical direct sum 
subalgebra of $\g \oplus \g$ and not the possible direct sum 
inside $\g$.)
Recall that the graph of $\pi_G^\sharp : T_{bw}^* G \ra T_{bw} G$
is given by Lemma 3.4. If $R_{bw}(x) \in \pi_G^\sharp
\left((T_{bw}(H U^-_w w))^0\right)$, then there exist
$n^+ \in \n^+ \cap \Ad_w(\n^-)$ and $n^- \in \n^-$ such that
$$
(x + n^+, x+ n^-) \in \Ad_{bw}(\Lie F) \subset \Ad_{bw}(\b^+ \oplus
\b^-).
$$
Comparing the first components gives 
$x \in \Ad_{bw}(\b^+) + \n^+ \cap \Ad_w(\n^-)$, and consequently
$$
\aligned
\pi_G^\sharp \left((T_{bw}(H U^-_w w))^0\right)
 &\subset R_{uw} \left(\Ad_{bw}(\b^+) + \n^+ \cap \Ad_w(\n^-)\right)  \\
 &= L_{bw}(\b^+) + R_{bw}(\n^+ \cap \Ad_w(\n^-)) =  \wt{E}_w(bw).
\endaligned
$$
Now the third statement in part (i) follows from Lemma 3.5.
\qed

\enddemo
\definition{3.7. Proof of Theorem 3.2} First we 
prove the Theorem in the case $J = \emptyset$.
For all $w \in W$, the tangent spaces to the 
fibers of $\nu_\emptyset$ inside $U^-_w y_w$ are 
$\eta_*(L_{bw}(\b^+))$ for $b \in H U^-_w$, and they 
are contained in $E_w$. Proposition 2.5 and part (2)
of Proposition 3.6 imply that
$$  
\pi_{G/H}|_{U_w^- y_w} - i_{w*}^\emptyset(\pi|_{B^- x_w})
\in  \wedge^2 E_w, \quad \text{for all} \; w \in W,
\tag 3.11
$$
which establishes the Theorem in the case $J=\emptyset$.

For the general case, in addition to Proposition 3.6 (ii), 
we need to prove that
$$
\pi_{G/H}|_{U_w^- y_w} - i_{w*}^J(\pi_J|_{B^- x_w^J})
\in  \wedge^2 E_w, \quad \text{for all} \; J \subset \Ga, \;
w \in W^J.
\tag 3.12
$$
Because $i^J_w = i^\emptyset_w \circ 
\left( \mu_J|_{B^- x_w} \right)^{-1}$ and 
$$
\mu_J|_{B^- x_w} : (B^- x_w, \pi|_{B^- x_w}) \ra 
(B^- x_w^J, \pi_J|_{B^- x_w^J})
$$
is a Poisson isomorphism for all $w \in W^J$
(recall (1.12) and Proposition 1.6), we get that
$$
i_{w*}^J(\pi_J|_{B^- x_w^J}) = 
i_{w*}^\emptyset(\pi|_{B^- x_w}),
\quad \text{for all} \; w \in W^J.
\tag 3.13
$$
Equations (3.11) and (3.13) imply (3.12), and this 
completes the proof of the Theorem.
\qed
\enddefinition

\proclaim{3.8. Lemma} The $H$-orbits of symplectic leaves
of $(G/H, \pi_{G/H})$ are exactly the projections 
$\tau(B^- w_1 B^- \cap B^+ w_2 B^+)$ of the 
double Bruhat cells of $G$ onto $G/H$, for $w_1,w_2\in W$.
\endproclaim

\demo{Sketch of the proof} The proof of the Lemma is 
analogous to the well known fact that the $H$-orbits
of symplectic leaves of $(G, \pi_G)$ are the double 
Bruhat cells of $G$.
The Drinfeld Lagrangian subalgebra of the base point $e H$ 
of $G/H$ is $\Lie (H_{\diag}(U^+ \times U^-))$,
see \cite{\BGY, Theorem 1.8}. The $H$-orbits of 
symplectic leaves of $(G/H, \pi_{G/H})$ are the inverse 
images of the $(B^+ \times B^-)$-orbits on 
$(G \times G)/H_{\diag}(U^+ \times U^-)$ under the map 
$$
G/H \ra (G\times G) / H_{\diag}(U^+ \times U^-), \quad
gH \mapsto (g,g) H_{\diag} (U^+ \times U^-).
\tag 3.14
$$ 
By the Bruhat Lemma, the $(B^+ \times B^-)$-orbits on 
$(G\times G) / H_{\diag}(U^+ \times U^-)$ are 
parametrized by $W \times W$ via
$$
W \times W \ni (w_1, w_2) \longmapsto 
(B^+ \times B^-) (w_2, w_1) H_{\diag} (U^+ \times U^-).
$$
Finally,
the inverse images of the above orbits under 
(3.14) are exactly the projections $\tau(B^- w_1 B^- \cap B^+ w_2 B^+)$. 
\qed
\enddemo

\definition{3.9. A second proof of Theorem 1.8} For
$w \in W^J$ and $w_1,w_2\in W$, the intersection
of $\tau(B^- w_1 B^- \cap B^+ w_2 B^+)$ with 
$i_w^J(B^- x_w^J) = U^-_w y_w$
is nonempty only if 
$w_1 = w$ and $w_2 \geq w$, because it lies 
inside $\tau(B^- w \cap B^- w_1 B^-)$ 
(thus $w_1=w$) and consequently inside 
$B^- w B^+ \cap B^+ w_2 B^+$ (thus $w_2 \geq w$ by 
\cite{\Deo, Corollary 1.2}).
If $w_1= w$ and $w_2 \geq w$, then the intersection
of $\tau(B^- w_1 B^- \cap B^+ w_2 B^+)$ with 
the image of $i_w^J$ is $\tau(U_{\dw, w_2})$, 
cf\. (0.6), and is a nonempty irreducible subvariety
of $G/H$ by \cite{\BGY, Theorem 2.3}.
Theorem 3.2, Lemma 3.8, and the argument 
of the proof of \cite{\BGY, Theorem 1.8} 
imply that for each $w_2 \in W$ with $w_2 \geq w$, the set
$\tau(\U_{\dw, w_2})$ is a (single) 
$H$-orbit of symplectic leaves of
$(i_w^J(B^- x_w^J), i_{w*}^J(\pi_J|_{B^- x_w^J}))$.
Thus, the $H$-orbits of symplectic leaves 
of $(B^- x_w^J, \pi_J|_{B^- x_w^J})$ are 
$$
(i_w^J)^{-1} \circ \tau (\U_{\dw, w_2}) = 
\U_{\dw, w_2} . P^+_J = \SS_{w, w_2}^J, \quad
w_2 \in W,\; w_2 \geq w.   
$$
This completes the second proof of the Theorem.
\qed 
\enddefinition

Let us note that although in the second proof of Theorem 1.8
we still used the theory of Poisson homogeneous spaces
in Lemma 3.8, we avoided the combinatorics 
arguments from the first proof. The latter were replaced
by the geometric construction of weak splittings of Poisson
submersions. In addition, 
in a subsequent publication we will demonstrate 
that those geometric arguments can be extended to the 
quantum situation, while we are not aware of any quantum 
version of the dressing orbit method that provides 
a classification of $H$-invariant prime ideals of some class of 
associative algebras.


\head 4. Hermitian symmetric spaces of 
compact type.
\endhead

\definition{4.1}
In this Section we investigate the Poisson 
structure $\pi_J$ on $G/P_J^+$ for the case 
when the unipotent radical $U_J^+$ of $P_J^+$
is abelian and $G$ is a simple algebraic 
group.
The flag varieties $G/P_J^+$ of this type 
exhaust all irreducible Hermitian symmetric 
spaces of compact type, \cite{\Wone}. 
We show that in this case all $L_J$-orbits on $G/P_J^+$ are
complete Poisson subvarieties with respect to 
the Poisson structure $\pi_J$. 
We then use the results of 
Richardson, R\"ohrle, and Steinberg \cite{\RRS}
on special representatives for 
the $L_J$-orbits on $G/P_J^+$. 
We prove that the 
Poisson structure $\pi_J$ vanishes at all
such special base points, and as a result of this
$(G/P_J^+, \pi_J)$ stratifies into complete 
Poisson subvarieties each of which is a 
quotient of $L_J$, equipped with the 
standard Poisson structure. 

First, set
$$
r_J = \sum_{\al \in \De_J^+} e_\al \wedge f_\al \quad
\text{and} \quad 
\pi_{L_J}= L(r_J)-R(r_J),
\tag 4.1
$$
cf\. (0.1) and (0.2).
These are respectively the standard $r$-matrix
and the standard Poisson structure on the reductive
group $L_J$. It is well known that $(L_J, \pi_J)$
is a Poisson algebraic subgroup of $(G, \pi_G)$.  
Set also
$$
\check{r}_J = \sum_{\al \in \De^+ \backslash \De_J^+} 
e_\al \wedge f_\al.
\tag 4.2
$$
Observe that $r_\g = r_J + \check{r}_J$.
\enddefinition

\proclaim{4.2. Proposition} Assume that $J \subset \Ga$ 
is such that $U_J^+$ is abelian. Then the following 
properties hold.

{\rm(i)} The Poisson structure $\pi_J$ on $G/P_J^+$ is 
given by 
$$
\pi_J= - \chi(r_J).
\tag 4.3
$$  
In particular, all $L_J$-orbits on the flag variety $G/P_J^+$ 
are complete Poisson subvarieties of $(G/P_J^+, \pi_J)$.

{\rm(ii)} Under the identification 
$\Psi_J : \n_J^- @>{\cong}>> B^- .P_J^+ \subset G/P_J^+,$ 
$\Psi_J(x) = \exp(x) P_J^+$,
of $L_J$-spaces {\rm(}where $L_J$ acts on the first term 
by the adjoint action{\rm)}, the restriction of 
the Poisson structure $\pi_J$ 
to $B^- . P_J^+$ corresponds to 
$-\chi(r_J) \in \Ga(\n_J, \wedge^2 T \n_J)$.
Here $\chi : \wedge^2 \l_J \ra \Ga(\n_J^-, \wedge^2 T \n_J^-)$
is derived from the adjoint action of $L_J$ on $\n_J^-$.
\endproclaim
 
\demo{Proof} (i) For all $u \in U_J^-$,
$$
L_u(\check{r}_J) - R_u(\check{r}_J)= 
R_u \left( \Ad_u \check{r}_J - \check{r}_J \right)=0
$$
because $U^-_J$ is abelian and 
$\check{r}^-_J \in \n^+_J \wedge \n^-_J$. Thus, 
$\pi_G(u) = \chi^R(r_J) - \chi^L(r_J)$ for 
$u \in U^-_J$. Since 
$\eta_{J*} \chi^R(r_J)=0$, we have $\pi_J= -\chi(r_J)$ on
the open $B^-$-orbit on $G/P_J^+$. But both bivector fields 
in (4.3) are 
algebraic, so they coincide.

The second statement in (i) directly follows from (4.3).
It is also a consequence of Theorem 1.1(b) in \cite{\RRS},
stating that each $L_J$-orbit on $G/P_J^+$ is the 
intersection of a $P_J^+$-orbit and a $P_J^-$-orbit.
(The latter are complete Poisson subvarieties of 
$(G/P_J^+, \pi_J)$ as was shown in \S 1.6.)

Part (ii) follows from part (i) by noting that 
$B^-. P^+_J = U^-_J. P^+_J$
is $L_J$ stable and $\Psi_J$ intertwines the actions 
of $L_J$ on $\n^-_J$ and $B^-.P^+_J$.  
\qed
\enddemo

\definition{4.3} The $L_J$-orbits on $G/P_J^+$ were classified
by Richardson, R\"ohrle, and Steinberg \cite{\RRS}, and previously
the $L_J$ orbits on the unipotent radical 
$U_J^-\cong B^-. P_J^+ \subset G/P_J^+$ were 
treated by Wolf \cite{\Wtwo} and Muller, Rubenthaller, 
and Schiffmann \cite{\MRS}. 
We recall the parametrization of the $L_J$-orbits on $G/P_J^+$
from \cite{\RRS, Theorem 1.2}. Fix a maximal
sequence $(\be_1, \ldots, \be_k)$ of long roots
in $\De^+\backslash \De^+_J$ which are mutually orthogonal.
Denote by $u_{-\be_i}$ a nontrivial element in 
the one-parameter unipotent subgroup of $G$ corresponding 
to the root $\be_i.$ Denote by $w_{\be_i} \in W$ 
the reflection corresponding to $\be_i.$

Note that $\{w_{\be_i}\}_{i=1}^k$ mutually commute because 
$\{\be_i\}_{i=1}^k$ are mutually orthogonal. 
(The elements
$\{u_{-\be_i}\}_{i=1}^k$ also mutually commute
because $U^-_J$ is abelian.)
Let $\dw_{\be_i}$ be a representative 
of $w_{\be_i}$ in the normalizer of the 
maximal torus $H$ of $G$. Finally, set
$$
x_{st} =\prod_{i=1}^t \dw_{\be_i} 
\prod_{j=t+1}^s u_{-\be_j} 
\quad \text{for} \quad
0 \leq t \leq s \leq k. 
$$
\enddefinition

\proclaim{Theorem} {\rm(Richardson--R\"ohrle--Steinberg \cite{\RRS})}
If the parabolic subgroup $P_J^+$ of $G$ has abelian unipotent
radical, then $\{ x_{st} P_J^+ \mid 0 \leq t \leq s \leq k\}$
is a system of representatives for the $L_J$-orbits 
on $G/P_J^+$.
\endproclaim

\definition{4.4} Let us recall that $P_J^+$
has an abelian unipotent radical only if 
it is a maximal parabolic subgroup of 
$G$. In addition, if 
$J = \Ga \backslash \{\ \al' \}$, then 
$P_J^+$ has abelian unipotent radical if 
and only the $\al'$ height of 
the longest root $\theta$ of $\g$ is equal
to 1. (We fix this root $\al'$ for the remainder
of this Section.) If this condition is satisfied,
then the $\al'$-height of any root 
$\ga \in \De$ is equal to $0, \pm 1$. 
Moreover
$n_{\al'}(\ga)= 1$, $0$, or $-1$,
if $\ga$ is a root of $\n_J^+$,
$\l_J$, or $\n_J^-$, respectively
(i.e., 
$\ga \in \De^+\backslash \De^+_J$, 
$\De_J$, or $-(\De^+\backslash \De^+_J)$).
\enddefinition

\definition{4.5} We will work with a 
special maximal set $(\be_1, \ldots, \be_k)$
of mutually orthogonal long roots in
$\De^+\backslash \De^+_J$. 
We proceed analogously to the proof of 
\cite{\RRS, Proposition 2.8}, defining
inductively $\be_i$ and subsets $\Ga_i$
of the set $\Ga$ of positive simple roots of $\g$.
Let $\be_1= \theta$ be the highest root of $\g$ 
and $\Ga_1=\Ga$. Assume that for 
some $i\leq k$, we have already defined $\be_i$ and $\Ga_i$. Let
$\wt{\Ga}_{i+1}$ be the set of  all roots in $\Ga_i$ that are orthogonal
to $\be_i.$ (Since $\be_i$ is dominant 
in the root system defined 
by $\Ga_i$, all roots in it that are orthogonal to 
$\be_i$ are combinations of simple roots
in $\wt{\Ga}_{i+1}$.)
If $k>1$, then $\wt{\Ga}_{i+1}$ contains 
$\al'$. Denote by $\Ga_{i+1}$ the connected 
component of $\wt{\Ga}_{i+1}$ containing $\al'$.
(Here we identify $\Ga$ with the Dynkin graph 
of $\g$ and view $\wt{\Ga}_{i+1}$ as a subgraph 
of it.) Finally, set $\be_{i+1}$ to be the 
highest root of the root system 
defined by $\Ga_{i+1}$.

This sequence has the properties that 
$$
\supp \be_j \subseteq \supp \be_i = \Ga_i 
\quad \text{for} \quad j \geq i
\tag 4.4
$$ 
(cf\. \S 0.7)
and that $\be_i$ is the longest root 
in the root system defined by $\Ga_i$.

For simplicity of the exposition, we set 
$$
u_{-\be_i} = \exp( f_{\be_i}),
\tag 4.5
$$
but we note that all proofs work for general $u_{-\be_i}$.

We will also use special representatives
for the reflections $w_{\be_i} \in W$ in the 
normalizer of the maximal torus $H$ of $G$
(and will thus omit the dot on top of them).
Set
$$
w_{\be_i} = \exp(e_{\be_i}) 
\exp(-f_{\be_i}) \exp(e_{\be_i}).
\tag 4.6
$$
This normalization is not necessary for the proof of the main
result in Theorem 4.6, but it simplifies the exposition.
Since $\be_i$ and $\beta_j$ are orthogonal for $i \neq j$ and 
$\be_i+\be_j$ is not a root of $\g$, we have that 
$w_{\be_i}(\be_j)= 0$ and  $[f_{\be_i}, f_{\be_j}] =0$, 
and consequently 
$$
\alignat2
[e_{\be_i}, f_{\be_j}] &= 0, 
 &\quad\qquad \Ad_{w_{\be_i}}(f_{\be_j}) &= f_{\be_j},  \\
\Ad_{w_{\be_i}}(u_{-\be_j}) &= u_{-\be_j},
 &\quad\qquad \Ad_{w_{\be_i}}(e_{\be_j}) &= e_{\be_j}
\endalignat
$$
for the special representatives (4.6) of $w_{\be_i}$.

The first main result in this sections is contained in 
the following Theorem.
\enddefinition

\proclaim{4.6. Theorem} If $P_J^+$ is a parabolic
subgroup with abelian
unipotent radical in the complex simple algebraic group $G$, then the
Poisson structure $\pi_J$ vanishes at the base points $x_{st} P_J^+$ of
the $L_J$-orbits  on $G/P_J^+$ $(0 \leq t \leq s \leq k)$. Therefore, all
$L_J$-orbits on $G/P_J^+$ are complete Poisson 
subvarieties of $(G/P_J^+, \pi_J)$ and are 
quotients of $(L_J, \pi_{L_J})$, cf\. 
\cite{\BGY, Theorem 1.8}.
\endproclaim 

Even in the case of Grassmannians, Theorem 4.6 contains
new information, compared to Part I \cite{\BGY}, 
where we dealt with the 
open $B^-$-orbits on Grassmannians. 

Let us note that Theorem 4.6 is not valid for 
an arbitrary maximal set $(\be_1, \ldots, \be_k)$
of mutually orthogonal long roots of 
$\De^+\backslash \De^+_J$. For instance,
in the case of $A_4$ and $\al'=\epsilon_3-\epsilon_2,$ 
take $\be_1= \epsilon_3 - \epsilon_1$ and
$\be_2 = \epsilon_4 - \epsilon_2$. Then 
$\pi_J$ does not vanish at 
$u_{-\be_1} u_{-\be_2} P_J^+$. The easiest way to see this 
is to note that the Poisson structure $\pi_{2,2}$
on the matrix affine Poisson space $M_{2,2}$ 
(cf\. \cite{\BGY, eq\. (1.7)}) does not vanish at 
$\leftmatrix 0 &1 \\ 1& 0 \rightmatrix$ and then to use
the embedding from \cite{\BGY, Proposition 3.4}.  

For the proof of Theorem 4.6, we will need several Lemmas.

\proclaim{4.7. Lemma} Assume that $\be$ is a long root 
in $\De^+\backslash\De_J^+$ and set
$$
w_\be = \exp(e_\be) \exp(-f_\be) \exp(e_\be).
$$
Let $\al$ be any root of $\g$ and $y_\al\in \g$ a nontrivial 
vector in the corresponding root space. If $[f_\be, y_\al]\neq 0$, 
then $w_\be(\al)= \al -\be$ and 
$$
[f_\be, y_\al] = - \Ad_{w_\be}(y_\al).
$$
\endproclaim  

\demo{Proof} Firstly, there is no root $\ga$ of $\g$ such that 
$\ga-\be$ and $\ga- 2 \be$ are roots as well. 
For if this happens, 
then the $\al'$ height of $\ga$ needs to be equal to 1, 
i.e., $\ga \in \De^+ \backslash \De^+_J$, because the 
$\al'$ heights of all roots of $\g$ are 0 or $\pm1$. 
Consequently,
$\ga + \be$ would not be a root and 
$\langle \ga, \be\spcheck \rangle \geq 2$. This would 
be a contradiction,
since $|\langle \ga, \be\spcheck \rangle|\leq 1$
because $\be$ is a long root. From this we get that 
the only roots of $\g$ of the form $\al + i \be$
for $i \in \ZZ$ are $\al$ and $\al -\be$,
thus $w_\be(\al)= \al -\be$, 
$\ad_{e_\be}(x_\al)=0$,
and $\ad_{f_\be}^2(x_\al)=0$. Consequently, 
$\Span \{y_\al, [f_\be, y_\al]\}$, under the adjoint 
action, 
is isomorphic to the vector representation of 
the $\sll_2$ triple 
$\{e_\be, \be\spcheck=[e_\be, f_\be], f_\be\}$.
By a standard computation, 
in this basis $w_\be$ acts by 
$\leftmatrix 0 & 1 \\ -1 & 0 \rightmatrix$, so
$\Ad_{w_\be}(y_\al)= - [f_\be, y_\al]$. 
\qed\enddemo

\proclaim{4.8. Lemma} If $\al \in \De_J^+$ is such that 
$$
[f_{\be_i}, e_\al] \neq 0 \quad \text{and}
\quad
[f_{\be_j}, f_\al] \neq 0
$$
for some $i \neq j \leq k$, 
then $\be_i - \al - \be_j \in \De_J^+$ and 
$$
[f_{\be_i},  e_\al] \wedge [f_{\be_j}, f_\al]+
[f_{\be_i}, e_{\be_i - \al - \be_j}] \wedge
[f_{\be_j}, f_{\be_i - \al - \be_j}]=0.
\tag 4.7
$$
\endproclaim

\demo{Proof} It follows from Lemma 4.7 that 
$$
w_{\be_i}(\al)=\al -\be_i \quad 
\text{and} \quad
w_{\be_j}(-\al) = -\al -\be_j.
$$
Since $\be_i$ and $\be_j$ are orthogonal,
$w_{\be_j}(\be_i)=0$ and 
$$
w_{\be_j} w_{\be_i}(\al)= w_{\be_j}(\al -\be_i)=
\al + \be_j - \be_i.
$$ 
Therefore $\be_i - \al - \be_j$ is a root of $\g$.
Because $\al - \be_i$ is a root of $\al'$ height 
equal to $-1$, it belongs to $-(\De^+\backslash \De_J^+)$.
As a consequence of this, 
$\supp \al \subset \supp \be_i = \Ga_i$, recall \S 4.5.
If $i < j$, then $\supp \al \subset \Ga_i \subset \Ga_j$,
and the fact that $\al + \be_j$ is a root would 
contradict  the property that $\be_j$ is the highest 
root of the root system of $\Ga_j$, cf\. \S 4.5. Thus,
$i>j$ and $\al+\be_j$ is a positive root
in the root system of $\Ga_i$. Since $\be_i$ is 
the highest root in this root system, $\be_i -\al - \be_j$ 
needs to be a positive root.

Because $\Ad_{w_{\be_j} w_{\be_i}}$ preserves the 
bilinear form $\langle.,.\rangle$ on $\g$, 
$$
\Ad_{w_{\be_j} w_{\be_i}}(e_\al) = a f_{\be_i-\al-\be_j}
\quad \text{and} \quad
\Ad_{w_{\be_j} w_{\be_i}}(f_\al) = a^{-1} e_{\be_i-\al-\be_j} 
\tag 4.8
$$
for some $a \in \CCx$. Equation (4.6) implies that
$\Ad_{w_{\be_i}}(f_{\be_i})=-e_{\be_i}$. Combining
the above facts and using that 
$\Ad_{w_{\be_j}}(e_{\be_i}) = e_{\be_i}$ (cf\. \S 4.5)
leads to
$$
\Ad_{w_{\be_j} w_{\be_i}}([f_{\be_i}, e_\al])=
 - a [e_{\be_i}, f_{\be_i -\al -\be_j}].
$$ 
Thus, $[e_{\be_i}, f_{\be_i -\al -\be_j}] \neq 0$, 
and consequently
$[f_{\be_i}, e_{\be_i -\al -\be_j}] \neq 0.$ 
Applying twice Lemma 4.7 and (4.8), we get 
that
$$
\aligned
[f_{\be_i}, e_{\be_i -\al -\be_j}] &=
- \Ad_{w_{\be_i}}(e_{\be_i -\al -\be_j})=
- a \Ad_{w_{\be_i}^2 w_{\be_j}}(f_\al) \\
&=a \Ad_{w_{\be_j}}(f_\al)= - a [f_{\be_j},f_\al].
\endaligned
$$ 
In the third equality, we used that $\Ad_{w_{\be_j}}(f_\al)$ is the 
lowest weight vector for the vector representation of 
the $\sll_2$ triple
$\{e_{\be_i}, \be_i\spcheck=
[e_{\be_i}, f_{\be_i}], f_{\be_i} \}$ 
(under the adjoint action) and in this representation
$w_{\be_i}^2$ acts by $-\id$, recall (4.6). Analogously,
$$
\aligned
[f_{\be_j}, f_{\be_i -\al -\be_j}] &=
- \Ad_{w_{\be_j}}(f_{\be_i -\al -\be_j})=
- a^{-1} \Ad_{w_{\be_i} w_{\be_j}^2}(e_\al) \\
&= a^{-1} \Ad_{w_{\be_i}}(e_\al)= - a^{-1} [f_{\be_i},e_\al].
\endaligned
$$
Hence,
$$
\aligned
[f_{\be_i}, e_{\be_i -\al -\be_j}] \wedge 
[f_{\be_j}, f_{\be_i -\al -\be_j}] &=
(- a [f_{\be_j},f_\al]) \wedge (- a^{-1} [f_{\be_i},e_\al])
\\
&= - [f_{\be_i},e_\al] \wedge [f_{\be_j},f_\al].
\qquad \qed
\endaligned
$$
\enddemo

\proclaim{4.9. Lemma} For all $0 \leq t \leq s \leq k$,
$$
\Ad_{u_{-\be_{t+1}} \ldots u_{-\be_s}}^{-1} (r_J)=
r_J - \sum_{j=t+1}^s \sum_{\al \in \De_J^+} 
\left( [f_{\be_j}, e_\al] \wedge f_\al +
e_\al \wedge [f_{\be_j}, f_\al] \right).
$$
\endproclaim

\demo{Proof} For arbitrary $0 \leq i \leq k$ and $0 \leq j \leq k$ and 
$\al \in \De_J$, the sum $\al + \be_i + \be_j$ is not a root of $\g$
since its 
$\al'$ height is $-2$. Taking into account that 
$u_{-\be_j} = \exp(f_{\be_j})$, we get
$$
\multline
\Ad_{u_{-\be_{t+1}} \ldots u_{-\be_s}}^{-1} (r_J) = \\
\sum_{\al \in \De_J^+}
\left( e_\al - [f_{\be_{t+1}}, e_\al] - \ldots - 
[f_{\be_s}, e_\al] \right)
\wedge
\left( f_\al - [f_{\be_{t+1}}, f_\al] - \ldots -
[f_{\be_s}, f_\al] \right).
\endmultline
$$
The Lemma will follow if we show that for all 
$s \leq i, j \leq k$,
$$
\sum_{\al \in \De_J^+} [f_{\be_i}, e_\al]
\wedge [f_{\be_j}, f_\al]=0.
$$
This is a corollary of Lemma 4.8, which implies that 
all $\al \in \De^+_J$ for which $[f_{\be_i}, e_\al] \neq 0$
and $[f_{\be_j}, f_\al] \neq 0$ can be grouped in pairs such 
the sum of the corresponding expressions 
$[f_{\be_i}, e_\al] \wedge [f_{\be_j}, f_\al]$ will be equal
to $0$.
\qed\enddemo

In the setting of \S 4.3, for each 
$0 \leq t \leq k$ set
$$
w_t= \prod_{i=1}^t w_i.
\tag 4.9
$$
The same notation will be used for the representative of
$w_t$ in the normalizer of $H$ in $G$ which 
is the product of the representatives (4.6).

\definition{4.10. Proof of Theorem 4.6} We will 
prove that
$$
\Ad_{w_t}^{-1} \left( [f_{\be_j}, e_\al]) \wedge f_\al \right)
\in \p_J^+ \wedge \g,
\quad \text{for} \; \al \in \De^+_J,\; 0 \leq t < j \leq k.
\tag 4.10
$$
Analogously, one shows that
$$
\Ad_{w_t}^{-1} \left( e_\al \wedge [f_{\be_j}, f_\al] \right)
\in \p_J^+ \wedge \g,
\quad \text{for} \; \al \in \De^+_J,\; 0 \leq t < j \leq k.
\tag 4.11
$$
It is clear that 
$$
\Ad_{w_t}^{-1} \left( e_\al \wedge f_\al \right) 
\in \n^- \wedge \n^+,
\quad \text{for} \; \al \in \De^+_J,\; 0 \leq t \leq k.
\tag 4.12
$$
Lemma 4.9, (4.10)--(4.12), and the commutativity of 
$w_t$ and $u_{-\be_{t+1}} \ldots u_{-\be_s}$ (cf\. \S 
4.5) imply that
$$
\Ad_{x_{st}} (r_J) \in \p_J^+ \wedge \g, \quad
\text{for} \; 0 \leq t \leq s \leq k,
$$
which is equivalent to the vanishing of $\pi_J$
at $x_{st} P_J^+$.

Thus, we are left with showing (4.10). We will make 
use of the following fact \cite{\RRS, Lemma 2.10 (b)}:
\roster
\item"(*)"
 {\it{For $\ga \in -(\De^+ \backslash \De^+_J)$, 
the set of all $\be_i$, $0 \leq i \leq t$, not orthogonal
to $\ga$ has cardinality $0$, $1$, or $2$, and accordingly
$w_t(\ga) \in -(\De^+ \backslash \De^+_J)$, 
$\De_J$, or $\De^+ \backslash \De^+_J$.}} 
\endroster

If $w_t^{-1}(-\al) \notin -(\De^+ \backslash \De^+_J)$,
then $\Ad_{w_t}^{-1}(f_\al) \in \p_J^+$ and we are done.
If $w_t(\al)^{-1} \in -(\De^+ \backslash \De^+_J)$, then 
applying (*) for $\ga =w_t^{-1}(-\al)$ one gets that there
exists $0 \leq i \leq t$ such that 
$w_t^{-1}(-\al)= w_{\be_i}(-\al) = -\al -\be_i$. 
This means that $\langle \be_i\spcheck, \al \rangle= -1$,
and consequently 
$\langle \be_i\spcheck, \al - \be_j\rangle= -1$
($\langle \be_i, \be_j \rangle = 0$ since $j>t \geq i$).
If $[f_{\be_j}, e_\al]=0$, then we are done. If 
$[f_{\be_j}, e_\al]\ne 0$, 
then using (*) again, this time for $\ga = \al -\be_j$, we get 
$w_t^{-1}(\al - \be_j) \notin -(\De^+ \backslash \De^+_J)$; 
therefore $\Ad_{w_t}^{-1} [f_{\be_j}, e_\al]) \in \p_J^+$. 
This establishes (4.10) and thus the Theorem.
\qed
\enddefinition

\definition{4.11} In the last part of this Section 
(\S 4.11--4.13), for a parabolic subgroup $P_J^+$ with 
abelian unipotent radical, we characterize the symplectic 
leaves of $(G, \pi_J^+)$ within each $L_J$-orbit. 
First we recall some facts on minimal length representatives
in Weyl groups. For two subset $I$ and $J$ of the 
set of positive simple roots $\Ga$ of $G$ denote
by ${}^I W_{\min}$ and ${}^I W^J_{\min}$ the set of (unique) 
minimal length representatives of the cosets in
$W_I \backslash W$ and $W_I \backslash W /W_J$.
Recall the following standard facts, see 
e.g. \cite{\Car, Lemma 4.3}.
\enddefinition

\proclaim{Lemma} For arbitrary $I, J \subset \Ga$
the following hold.

{\rm(i)} Every element of $W^J_{\min}$ can be uniquely 
represented as a product $vw$ for some 
$w \in {}^I W^J_{\min}$ 
and $v \in W_I \cap W^{I \cap w(J)}_{\min}$.

{\rm(ii)} Every element of $W$ can be uniquely
represented as a product $v_1 w v_2$ for some
$w \in {}^I W^J_{\min}$, $v_2 \in W_J$ and
$v_1 \in W_I \cap W^{I \cap w(J)}_{\min}$.
\endproclaim

\definition{4.12} We will need the following results 
from \cite{\RRS}, recall (4.9).
\enddefinition

\proclaim{Theorem} {\rm(Richardson--R\"ohrle--Steinberg)}
If $P_J^+ \subset G$ has abelian unipotent radical,
then the following hold.

{\rm(i)} For any given $P_J^+$-orbit and any given $P_J^-$-orbit 
on $G/P_J^+$, the two orbits are either disjoint or else intersect 
in a single $L_J$-orbit. 

{\rm(ii)} The set $\{w_s\}_{s=0}^k$ coincides with 
${}^J W^J_{\min}$ and thus 
$G/P_J^+ = 
\bigsqcup \Sb 0\leq s\leq k \endSb P_J^+. w_s P_J^+ =
\bigsqcup \Sb 0\leq s\leq k \endSb P_J^-. w_s P_J^+$.

{\rm(iii)} For all $0 \leq t \leq s \leq k$, we have
$$  
L_J . x_{st} P_J^+ = P_J^- . w_t P_J^+ \cap
P_J^+ . w_s P_J^+,
$$
cf\. {\rm (4.9)}.
\endproclaim 

Part (i) is \cite{\RRS, Theorem 1.1 (b)}, part (ii) is
\cite{\RRS, Proposition 2.11}, and part (iii) is 
\cite{\RRS, Lemma 3.5 (d)}. We will illustrate the inclusion
$$
L_J . x_{st} P_J^+ \subset P_J^- . w_t P_J^+ \cap
P_J^+ . w_s P_J^+.
$$ 
Since $\dw_{\be_i}$ and $u_{-\be_j}$ commute for 
$i \neq j$ (cf\. \S 4.5),
$$
L_J . x_{st} P_J^+ = 
L_J \prod_{j=t+1}^s u_{-\be_j}. 
\prod_{i=1}^t \dw_i P_J^+ \subset P_J^- w_t P_J^+.
$$
Because $[e_{\be_i},f_{\be_j}]=0$ for
$i \neq j$ and $w_{\be_i}^2=1$,
$$
\aligned
L_J. x_{st} P_J^+ &= L_J \prod_{i=1}^t \exp(e_{\be_i})  
\left( \prod_{i=1}^t \exp(-e_{\be_i}) \exp(f_{\be_i})
\exp(-e_{\be_i}) \right) \prod_{j=t+1}^s w_{\be_j} P_J^+  \\
 &\subset P_J^+ . w_s P_J^+.
\endaligned
$$  

The second main result of this section is contained in
the following Theorem. 

\proclaim{4.13. Theorem} Assume that $P_J^+$ is 
a parabolic subgroup of $G$ with abelian unipotent
radical. For all $0 \leq s \leq t \leq k$, the $H$-orbits
of symplectic leaves of 
$(L_J. x_{st} P_J^+,
\pi_J|_{L_J. x_{st} P_J^+})$
are $\SS^J_{(v_1 w_t \wJ, v_2 w_s v_3)}$ {\rm(}cf\. {\rm(1.9))}
for unique 
$$
\xalignat3
v_3 &\in W_J,  &v_1 &\in W_J \cap W^{J \cap w_t(J)}_{\min},
&v_2 &\in W_J \cap W^{J \cap w_s(J)}_{\min}
\tag 4.13  \endxalignat  
$$
 such that
$v_1 w_t \wJ \leq v_2 w_s v_3$.
\endproclaim 

As an example, let us note that the $L_J$-orbits 
inside the open $B^-$-orbit on $G/P_J^+$ are the orbits 
$L_J .x_{0s} P_J^+$, for $0 \leq s \leq k$. 
Theorem 4.13 implies that the symplectic
leaves of $(L_J .x_{0s} P_J^+, \pi_J|_{L_J .x_{0s} P_J^+})$
are $\SS^J_{(\wJ, v_2 w_t v_3)}$ for some 
$v_3 \in W_J$ and
$v_2 \in W_J \cap W^{J \cap w_s(J)}_{\min}$ such that
$\wJ \leq v_2 w_s v_3$. 

\demo{Proof} Because 
$\U_{\dw_1, w_2} \subset B^- w_1 \cap B^+ w_2 B^+$, recall 
(0.6), we have
$$
\SS^J_{w_1, w_2} = \U_{\dw_1, w_2}. P_J^+ \subset
B^- x_{w_1}^J \cap B^+ x_{w_2}^J \subset 
P_J^- x_{w_1}^J \cap P_J^+ x_{w_2}^J.
$$
It follows from Lemma 4.11, Theorem 4.12 (ii), and 
(1.3) that every element of $\Om^J$ can be uniquely
represented in the form 
$(v_1 w_t \wJ, v_2 w_s v_3)$ for some 
$0 \leq s, t \leq k$ and $v_3$, $v_1$, $v_2$ as in (4.13) such that
$v_1 w_t \wJ \leq v_2 w_s v_3$. For such a 
pair, we have 
$$
\SS^J_{(v_1 w_t \wJ, v_2 w_s v_3)} \subset 
P_J^-. v_1 w_t \wJ P_J^+ \cap 
P_J^+. v_2 w_s v_3 P_J^+
= P_J^+. w_t P_J^+ \cap P_J^- . w_s P_J^+.
$$
Theorem 4.12 implies that the last 
intersection is nontrivial only if $t \leq s$, 
in which case it is equal to $L_J. x_{ts} P_J^+$. 
This completes the proof of the Theorem. 
\qed\enddemo


\head 5. The open $B^-$-orbits in compact Hermitian 
symmetric spaces
\endhead

In this section, we treat in detail the
restriction of the Poisson structure 
$\pi_J$ to the open $B^-$-orbit of $G/P_J^+$ 
in the case when $G$ is a complex simple group and $P_J^+$ 
is a parabolic subgroup with abelian unipotent radical. 
First, in \S 5.1, we obtain general formulas for 
$\pi_J|_{B^-. P_J^+}$.
Then, in \S 5.3--5.6, we use it to derive explicit formulas 
for all classical groups and show that all 
such Poisson structures are quasiclassical limits of 
interesting classes of quadratic algebras, of the type known as {\it
quantized coordinate rings\/} of classical varieties. (See, e.g.,
\cite{\G} for a general survey of quantized coordinate rings.) In
\S 5.7, we show that the
$E_6$ case gives rise  to a new quadratic Poisson structure on a 16 
dimensional affine space, related to a half-spin 
representation of $\so_{10}$.
 
\definition{5.1} For a general reductive Lie algebra $\g$
with fixed dual Borel subalgebras $\b^\pm$ as in 
the Introduction, consider an irreducible representation
$V_\g^\lambda$ of $\g$ with highest weight $\lambda$. Define 
the bivector field 
$$
\pi_\g^\lambda = \chi (r_{\g}) \in \Ga( (V_\g^\lambda)^*, 
\wedge^2 T (V_\g^\lambda)^*)
\tag 5.1
$$
(cf\. (0.1) and \S 0.2) where $\chi$ is derived from the 
action of $\g$ on $(V_\g^\lambda)^*$.
If $\pi_\g^\lambda$ is Poisson, 
then the corresponding Poisson bracket on the algebra of regular 
functions on $(V_\g^\lambda)^*$, identified 
with the symmetric algebra $S(V_\g^\lambda)$, is 
induced by 
$$
\{v_1, v_2\} = m( r_{\g} (v_1 \otimes v_2)), \quad v_1, v_2 \in 
V_\g^\lambda
\tag 5.2
$$
where 
$m : V_\g^\lambda \otimes V_\g^\lambda \ra S(V_\g^\lambda)$   
is the multiplication map, $m(v_1, v_2)= v_1 v_2$.

In the setting of \S 0.2, for $\al \in \De_J^+$ and 
$\be \in \De^+ \backslash \De^+_J$ set
$$
[e_\al, e_\be] = N_{\al, \be} e_{\al+\be}, \quad 
[f_\al, e_\be] = N_{-\al,\be} e_{-\al+\be}, \quad 
N_{\al,\be}, N_{-\al,\be} \in \CC.
\tag 5.3
$$

For a Lie algebra $\g$, we will denote by
$\g'=[\g, \g]$ its derived subalgebra. 

Finally, for arbitrary $J \subset \Ga$, 
denote by $\theta|_J$ the restriction of the 
highest root of $\g$ to (a dominant weight) of $\l_J$.
\enddefinition

\proclaim{Proposition} Assume that $G$ is a complex simple algebraic 
group and $P_J^+$ is a parabolic subgroup of $G$ with 
abelian unipotent radical. Then the following hold.

{\rm(1)} Under the adjoint action of $\l_J$, the space $\n_J^+$ is
an irreducible representation of $\l_J'$ with highest
weight $\theta|_J$. The restriction of the 
Poisson structure $\pi_J$ to $B^-. P_J^+ \subset G/P_J^+,$
identified with $\n^-_J$ by $x \mapsto \exp(x) P_J^+$, coincides
with $-\pi_{\l'_J}^{\theta|_J}$, cf\. {\rm(5.1)--(5.2)}. 

{\rm(2)} The Poisson structure on 
$\n_J^- = \{ 
\sum_{\be \in \De^+ \backslash \De^+_J} y_\be f_\be
\mid y_\be \in \CC\} \cong \AA^{|\De^+| - |\De^+_J|}$
from the first part is also is given by
$$
\{y_\be, y_\ga\} = \sum_{\al \in \De^+_J}
\bigl( -N_{\al, \be} N_{-\al,\ga} 
y_{\al+\be}  y_{-\al+\ga} +
N_{-\al, \be} N_{\al,\ga} 
y_{-\al+\be}  y_{\al+\ga} \bigr), 
\quad \be, \ga \in \De^+ \backslash \De^+_J.
$$
\endproclaim

\demo{Proof} The first statement of part (1) 
is well known. The simplest way to show it
here, is to observe that $L_J$ has finitely many 
orbits on $\n^-_J$ (acting by the adjoint action)
because it has finitely many orbits on $G/P_J^+$. 
This implies that $\n_J^-$ is an irreducible 
representation of $\l'_J$ (under the adjoint action)
and it must have highest weight $\theta|_J$. 

The second statement in part (i) now follows from 
the first one, Proposition 4.2(ii), and the 
definition (5.1).

Part (ii) is an immediate corollary of part (i), 
see (5.2). \qed
\enddemo

\definition{5.2}
Assume that the positive simple roots 
$\al_1, \ldots, \al_N$ of 
the simple algebraic group
$G$ are enumerated as in \cite{\Bour}. Recall that only 
maximal parabolic subgroups $P_J^+$ of $G$ can have an 
abelian unipotent radical. Moreover, those are exactly
the parabolic subgroups $P_{\Ga \backslash \{\al_m\}}^+$
for which the simple root $\al_m$ appears with
multiplicity 1 in the expansion of the highest root 
$\theta$ of $\g$ in terms of the positive simple roots,
i.e., the $\al_m$ height of $\theta$ is equal to 1.
In other words, $\al_m$ should be a cominuscule root.
This leads to the following 
choices for $\al_m$ according to the type of $G$:
$$
\alignat2
 &A_N &\qquad &\al_1, \ldots, \al_N,  \\
 &B_N &\qquad &\al_1,  \\
 &C_N &\qquad &\al_N,  \\
 &D_N &\qquad &\al_1, \al_{N-1}, \al_N,  \\
 &E_6 &\qquad &\al_1, \al_6,  \\
 &E_7 &\qquad &\al_7.
\endalignat
$$

For the other types ($E_8$, $F_4$, $G_2$), there are no 
parabolic subgroups with abelian unipotent radicals.
Below we will denote by $\omega_m$ the $m$-th 
fundamental weight of $G$.
\enddefinition

\definition{5.3. $A_N$ case} In this case, $\g = \sll_N$,
$$
\theta= \al_1 + \ldots + \al_N =
\omega_1 + \omega_N,
$$
and all simple roots $\al_m$
($m = 1, \ldots, N$) are cominuscule.
The derived subalgebra of the Levi factor 
of the parabolic subalgebra $\p^+_{\Ga \backslash \{\al_m\}}$ is 
$\l'_{\Ga \backslash \{\al_m\}} \cong \sll_m \oplus \sll_{N-m+1}$.
Under the adjoint action of $\l'_{\Ga \backslash \{\al_m\}}$, 
the nilradical 
$\n_{\Ga \backslash \{\al_m\}}^+$ is identified with 
$V_{\sll_m}^{\omega_1} \otimes 
\bigl( V_{\sll_{N-m+1}}^{\omega_1} \bigr)^*$
where $V_{\sll_l}^{\omega_1}$ denotes the vector representation 
of $\sll_l$. In Part I, 
we showed that, after applying the twist \cite{\BGY, eq\. (3.14)}, 
the induced Poisson structure on $\n_{\Ga \backslash \{\al_m\}}^+$ 
corresponds to the standard quadratic Poisson structure 
on the space of rectangular matrices $M_{N-m+1, m}$ (cf\. \cite{\BGY,
eq\. (1.7)}). It is the quasiclassical limit of the 
quadratic algebra of {\it{quantum matrices}}. 
\enddefinition

\definition{5.4. $B_{N+1}$ case} In this case, $\g \cong \so_{2N+1}$, 
$$
\theta= \al_1 + 2 \al_2 + \ldots + 2 \al_N= \omega_2,
$$
and the only root that is cominuscule is $\al_1$.  
Furthermore, $\l'_{\Ga\backslash \{\al_1\}} \cong \so_{2N-1}$,
and under the adjoint action the nilradical 
$\n_{\Ga \backslash \{\al_1\}}^+$ corresponds to 
the vector representation $V_{\so_{2N-1}}^{\omega_1}$ 
of $\so_{2N-1}$. If we identify 
$$
\multline
\n^-_{\Ga \backslash \{\al_1\}} \cong \bigl\{ 
\sum_{i=2}^N x_i (E_{N+i,1} - E_{N+1,i})
+\sum_{j=2}^N y_j (E_{j,1} - E_{N+1,N+j})
\\
+ z (E_{2N+1, 1} - E_{N+1, 2N+1})
\bigm| x_i, y_j, z \in \CC\bigr\}= \AA^{2N-1},
\endmultline
$$
then the induced Poisson structure on 
$\n^-_{\Ga \backslash \{\al_1\}}$ is given by
$$
\alignedat2
2\{x_i,x_j\} &= x_ix_j,  \qquad &&i<j                     \\
2\{y_i,y_j\} &= -y_iy_j,  \qquad &&i<j                     \\
2\{z,x_i\}   &= -zx_i,    \qquad &&\text{for all} \; j     \\
2\{z,y_j\}   &= zy_j,    \qquad &&\text{for all} \; j     \\
2\{x_i,y_j\} &= x_iy_j,   \qquad &&i \ne j                 \\
2\{x_i,y_i\} &= -2 \sum_{l > i} x_l y_l - z^2
                          \qquad &&\text{for all} \; i.
\endalignedat
$$
This is the quasiclassical limit of the {\it{odd dimensional
quantum Euclidean space}\/}, $O^{2N-1}_{q^{1/2}}(\CC)$, introduced by
Resh\-e\-tikh\-in, Takhtadzhyan, and Faddeev, \cite{\RTF, Definition
12}. (See \cite{\Mus, \S\S2.1, 2.2} for a simplified set of relations.)
\enddefinition

\definition{5.5. $C_{N+1}$ case} 
In this case, $\g \cong \sp_{2N}$ and
$$
\theta= 2\al_1 + \ldots + 2 \al_{N-1} + \al_N= 2 \omega_1.
$$
The only cominuscule root is $\al_N$, and  
$\l'_{\Ga\backslash \{\al_N\}} \cong \sll_N$.
Then under the adjoint action the nilradical 
$\n_{\Ga \backslash \{\al_N\}}^+$ corresponds to 
the second symmetric power 
$S^2(V_{\sll_N}^{\omega_1}) \cong V_{\sll_N}^{2 \omega_1}$
of the vector representation of $\sll_N$. If we
identify
$\n^-_{\Ga \backslash \{\al_N\}}$
with the space of symmetric matrices of 
size $N$ by
$$
\n^-_{\Ga \backslash \{\al_N\}} = 
\left\{ \left[ \smallmatrix 0 & 0 \\ Y & 0
\endsmallmatrix \right] \biggm| Y= (y_{ij})_{i,j=1}^N, \; 
y_{ij} = y_{ji} \right\},
$$
then the Poisson structure on 
$\n^-_{\Ga \backslash \{\al_N\}}$
is given by
$$
\aligned
\{y_{ij}, y_{lm}\} &= 
\bigl( \sign(m-i) + \sign(l-j) \bigr) y_{il} y_{jm}  \\
 &\qquad + \bigl( \sign(l-i) + \sign(m-j) \bigr) y_{im} y_{jl}.
\endaligned
$$
Interestingly, this Poisson structure is not, as one might expect after
seeing \S5.4, the quasiclassical limit of a quantum symplectic space
(see \cite{\RTF, Definition 14} or \cite{\Mus, \S1.1}). Instead, it
is the quasiclassical limit of the algebra of {\it quantum
symmetric matrices\/} introduced by Noumi in \cite{\Nou, Theorem 4.3,
Proposition 4.4, and comments following the proof} (with the parameters
$a_k$ all set equal to 1), and of the one given by Kamita \cite{\Kam, Theorem 0.2} (with $q$ and $q^{-1}$ interchanged).
\enddefinition

\definition{5.6. $D_{N}$ case} In this case, $\g =\so_{2N}$,
$$
\theta = \al_1 + 2 \al_2 + \ldots + 2 \al_{N-2} 
+ \al_{N-1} + \al_N = \omega_2,
$$
and the cominuscule roots are $\al_1$,
$\al_{N-1}$, and $\al_N$. (We will 
assume $N \geq 4$.) Below we consider 
those three cases.

{\bf{(a) Root $\al_1$.}} 
Here $\l'_{\Ga\backslash\{\al_1\}} \cong \so_{2(N-1)}$
and $\n^+_{\Ga\backslash\{\al_1\}}$, considered 
an $\l'_{\Ga\backslash\{\al_1\}}$-module under the adjoint 
action, is isomorphic to the vector representation
$V_{\so_{2(N-1)}}^{\omega_1}$
of $\so_{2(N-1)}$. We identify
$$
\aligned
\n^-_{\Ga \backslash \{\al_1\}} &= \bigl\{ 
\sum_{i=2}^N x_i (E_{N+i,1} - E_{N+1,i})
+\sum_{j=2}^N y_j (E_{j,1} - E_{N+1,N+j})
\bigm| x_i, y_j \in \CC\bigr\}  \\
 &= \AA^{2(N-1)}.
\endaligned
$$
Then the induced Poisson structure on 
$\n^-_{\Ga \backslash \{\al_1\}}$
is given by
$$
\alignedat2
2\{x_i,x_j\} &= x_ix_j,  \qquad &&i<j                     \\
2\{y_i,y_j\} &= -y_iy_j,  \qquad &&i<j                     \\
2\{x_i,y_j\} &= x_iy_j,   \qquad &&i \ne j                 \\
 \{x_i,y_i\} &= -\sum_{l > i} x_l y_l,
                          \qquad &&\text{for all} \; i.
\endalignedat
$$
It is the quasiclassical limit of the {\it{even dimensional
quantum Euclidean space}\/} $O^{2(N-1)}_{q^{1/2}}(\CC)$, see
\cite{\RTF,
\Mus}.

{\bf{(b) Root $\al_{N}$.}}
In this case,
$\l'_{\Ga\backslash\{\al_N\}} \cong \sll_{N}$
and under the adjoint action $\n^+_{\Ga\backslash\{\al_N\}}$
is isomorphic to 
$V_{\sll_N}^{\omega_2} \cong \wedge^2 V_{\sll_N}^{\omega_1}$ 
where $V_{\sll_N}^{\omega_1}$ is the vector representation
of $\sll_{N}$. Identify
$\n^-_{\Ga \backslash \{\al_N\}}$
with the space of skew-symmetric matrices of 
size $N$ by
$$
\n^-_{\Ga \backslash \{\al_N\}} = 
\left\{ \left[ \smallmatrix 0 & 0 \\ Y & 0
\endsmallmatrix \right] \biggm| Y= (y_{ij})_{i,j=1}^N, \; 
y_{ij} = - y_{ji} \right\}.
$$
Then the Poisson structure on 
$\n^-_{\Ga \backslash \{\al_N\}}$
is given by
$$
\aligned
2 \{y_{ij}, y_{lm}\} &= 
\bigl( \sign(l-j) + \sign(m-i) \bigr) y_{il} y_{jm}  \\
 &\qquad - \bigl( \sign(l-i) + \sign(m-j) \bigr) y_{im} y_{jl}.
\endaligned
$$
This is the quasiclassical limit of the algebra
of {\it{quantum antisymmetric matrices}\/} introduced 
by Strickland in \cite{\St, Section 1} (with $q$ replaced by $q^{1/2}$).

{\bf{(c) Root $\al_{N-1}$.}} One can lift the involutive 
automorphism of the Dynkin graph $D_N$ (preserving 
$\al_1, \ldots, \al_{N-2}$, 
and interchanging $\al_{N-1}$ and $\al_N$)
to an automorphism of $\so_{2N}$
that restricts to an isomorphism between 
$\l_{\Ga\backslash\{\al_{N-1}\}}$ and
$\l_{\Ga\backslash\{\al_N\}}$
(interchanging their $r$-matrices)
and their modules 
$\n^-_{\Ga \backslash \{\al_{N-1}\}}$ and 
$\n^-_{\Ga \backslash \{\al_{N}\}}$.
As a result, that 
automorphism of $\so_{2N}$ restricts to 
an isomorphism between the Poisson structures 
on $\n^-_{\Ga \backslash \{\al_{N-1}\}}$ and 
$\n^-_{\Ga \backslash \{\al_{N}\}}$.
\enddefinition

\definition{5.7. $E_6$ case} With this example, we show that
compact Hermitian symmetric spaces for exceptional groups
give rise to new interesting quadratic Poisson structures 
on affine spaces.

The highest root of the simple Lie algebra $\e_6$ of type $E_6$
is 
$$
\theta = \al_1 + 2 \al_2 + 2 \al_3 + 3 \al_4 + 2 \al_5 + \al_6
= \omega_2,
$$ 
and $\e_6$ has two cominuscule roots: $\al_1$ and $\al_6$.
Similarly to \S 5.6(c), one lifts the involutive automorphism 
of the Dynkin graph $E_6$ that interchanges $\al_1$ and $\al_6$
and fixes the other nodes to an automorphism of $\e_6$
that restricts to an isomorphism between
$\l'_{\Ga \backslash \{\al_1\}}$ and 
$\l'_{\Ga \backslash \{\al_6\}}$, and between 
their modules 
$\n^-_{\Ga \backslash \{\al_1\}}$ and 
$\n^-_{\Ga \backslash \{\al_6\}}$.
This linear map provides an isomorphism 
between the induced Poisson structures 
on $\n^-_{\Ga \backslash \{\al_1\}}$ and 
$\n^-_{\Ga \backslash \{\al_6\}}$.

In the case of the root $\al_1$, we have
$\l'_{\Ga\backslash\{\al_1\}} \cong \so_{10}$,
and as an $\l'_{\Ga\backslash\{\al_1\}}$-module 
$\n^+_{\Ga\backslash\{\al_1\}}$ is isomorphic to
one of the half-spin representations 
$V_{\so_{10}}^{\omega_5}$ of $\so_{10}$.
As a vector space, $V_{\so_{10}}^{\omega_5}$
is identified with
$$
V \oplus \wedge^3 V \oplus \wedge^5 V
$$
for a 5 dimensional vector space $V$, cf\. 
\cite{\FH, Section 20}. A basis 
$\{v_1, \ldots, v_5\}$ of $V$ gives rise 
to the basis 
$\{ v_{i_1} \wedge \dots \wedge v_{i_n} \mid
\; n\; \text{odd}, \; \; i_1 < \ldots <i_n \}$
of $V_{\so_{10}}^{\omega_5}$. We view it as
a set of coordinate functions
$$
\{ y_I \mid I \subset \{1, \ldots, 5\}, \; |I|\; \text{odd} \} 
$$
on $(V_{\so_{10}}^{\omega_5})^*$.
In terms of those coordinates, the induced quadratic Poisson 
structure on this 16 dimensional affine space is given by
$$
\multline
\{y_I, y_J \} =
\sum \Sb i \in I \backslash J \\ j \in J \backslash I \\ i \ne j \endSb
\sign(j-i) a_{i,j}^{I,J} y_{(I \backslash \{i\}) \cup \{j\}}
y_{(J \backslash \{j\}) \cup \{i\}}  \\
 -\frac{1}{2}
\sum \Sb \{i,j\} \subset I \backslash J \\ i < j \endSb
a_{i,j}^{I,J} y_{I \backslash \{i,j\}}
y_{J \cup \{i,j\}}  +  \frac{1}{2}
\sum \Sb \{i,j\} \subset J \backslash I \\ i < j \endSb
a_{i,j}^{I,J} y_{I \cup \{i,j\}}
y_{J \backslash \{i,j\}}.
\endmultline
$$
Here, for two subsets $I$ and $J$, and two elements $i$ and $j$ 
of $\{1, \ldots, 5\}$, we set
$$
I(i,j)= \{ l \in I \mid i<l <j \; \; \text{or} \; \; 
j < l < i \}
$$
and 
$$
a_{i,j}^{I,J} = (-1)^{|I(i,j)| + |J(i,j)|}.
$$
For example, for $i < j$,
$$
\{ y_{\{i\}}, y_{\{j\}} \} = y_{\{i\}} y_{\{j\}}.
$$
\enddefinition


\head Acknowledgements\endhead

We thank T. H. Lenagan, J.-H. Lu, K. Rietsch, 
R. Steinberg, and J. Wolf for helpful 
discussions and correspondence.


\enddocument